\documentclass[preprint,11pt]{elsarticle}

\usepackage{geometry}
\usepackage{amssymb}   %latex
\usepackage{amsmath}
\usepackage{float}
\usepackage[colorlinks,linkcolor=red]{hyperref}

\usepackage[ruled,lined,boxed,commentsnumbered]{algorithm2e}

\usepackage{graphicx}
\usepackage{subfigure}
\usepackage{color}
\usepackage{hyperref}
\usepackage{bm}

\usepackage{cases}

\def\dfrac{\displaystyle\frac}
\def\dint{\displaystyle\int}
\def\dsum{\displaystyle\sum}

\def\ud{\textrm{d}}

\journal{Journal}

\begin{document}
\begin{frontmatter}

\title{Time-marching multi-level variational multiscale tensor decomposition algorithm for heat conduction with moving heat source}

\author[PKU]{Xinyi Guan}
\author[PKU]{Jiayi Hu}
\author[MECH]{Lei Zhang}
\author[PKU]{Shaoqiang Tang\corref{cor1}}
\author[NU]{Wing Kam Liu\corref{cor1}}
\address[PKU]{HEDPS and LTCS, School of Mechanics and Engineering Science, Peking University, Beijing 100871, China}
\address[MECH]{School of Engineering Science, University of Chinese Academy of Sciences, Beijing 100049, China}%Institute of Mechanics, Chinese Academy of Sciences, Beijing 100080, China
\address[NU]{Co-founder of HiDeNN-AI, LLC; Department of Mechanical Engineering, Northwestern University, Evanston, IL 60208-3111, USA}
\cortext[cor1]{Author for correspondence, maotang@pku.edu.cn, w-liu@northwestern.edu}

\begin{abstract}

In this paper, we propose a time-marching multi-level \underline{V}ariational \underline{M}ulti\underline{s}cale-\underline{T}ensor \underline{D}ecomposition (VMS-TD) algorithm to solve the heat equation with a moving heat source model that arises from additive manufacturing. First, we take a second-order centered difference for time semi-discretization. The temperature field is decomposed according to multiple space resolution levels, each represented by the TD method. Then we adopt the VMS formulation \cite{Hughes1998} for the resulting elliptic problem to obtain a Galerkin weak form, and design VMS-TD algorithm to solve it. Furthermore, to comply with the TD solution scheme, special inter-scale data transfers are made at the scale interface and moving fine-scale subdomains. Numerical results demonstrate that the multi-level VMS-TD algorithm is much more efficient than the fully resolved TD algorithm, let alone traditional direct numerical simulation methods such as finite difference or finite element analysis. Compared with the well-known multi-grid methods or more recent GO-MELT framework\cite{Wagner2024}, the three-level VMS-TD uses much smaller degrees of freedom to reach accurate results. A multi-time-scale extension of VMS-TD algorithm is also proposed.
	
\end{abstract}
\begin{keyword}
%% keywords here, in the form: keyword \sep keyword
Multi-level variational multi-scale algorithm, Tensor decomposition, Model reduction, Heat conduction, Additive manufacturing
%% MSC codes here, in the form: \MSC code \sep code
%% or \MSC[2008] code \sep code (2000 is the default)
\end{keyword}

\end{frontmatter}

%%
%% Start line numbering here if you want
%%
% \linenumbers

%% main text

% \textbf{Nomenclature}

% 1. Continuous functions are denoted by lowercase letters, such as $u$, $f$;

% 2. Discretized functions (arrays, matrices) are denoted by uppercase letters, such as $U$, $F$;

% 3. Mappings (or certain operations in numerical implementation) are denoted by calligraphic letters, such as $\mathcal{E}$, $\mathcal{I}$;

% 4. Coarse and fine scale array indices are represented by uppercase and lowercase letters, respectively, such as $\overline{U_{\bm{J}}}$, $U_{\bm{j}}$;

% 5. Bold indices represent vectors, such as $U_{\bm{j}} = U_{j_x j_y}$.

\section{Introduction}

Additive manufacturing, also known as 3D printing, is a transformative technology that constructs objects layer by layer from digital models, allowing unparalleled design freedom and material efficiency. Originally developed for rapid prototyping in the 1980s, additive manufacturing has evolved into a suite of techniques with applications ranging from aerospace and automotive to biomedical and consumer goods industries \cite{Gani2025,Latiyan2025,Wang2025,Zhang2025_AM}. One popular type of additive manufacturing is laser powder bed fusion (LPBF)\cite{Ma2024}, which uses a high-powered laser to selectively melt and fuse layers of metal powder to create complex, high-precision components directly from digital models. In LPBF, there is an interplay among three physical objects, namely, the laser described by electrodynamics, the melt pool described by fluid mechanics, and the powder bed described by heat conduction. Here we only focus on heat conduction in the powder bed. Typical laser beam size ranges from 50 $\mu$m to 200 $\mu$m, depending on the machine and application \cite{Faegh2024}, while the powder bed size also varies, depending on whether it is designed for small components (e.g., dental implants) or large industrial parts. Typical commercial LPBF machine size is about 250mm $\times$ 250mm $\times$ 300mm \cite{Fardan2025}. In this case, the diameter of powder bed is $10^3$ to $10^4$ times larger than the laser beam. It causes immense computational cost if one performs direct numerical simulations by brute force for the heat conduction in the powder bed, and multi-scale method becomes indispensable for effective or even real-time simulations.

Almost half century ago, multi-grid/multi-level methods were proposed and developed in solving effectively discrete algebraic systems for boundary value problems that span over different scales \cite{Brandt1977,Hackbusch1985}.  Then over the past three decades, many multiscale frameworks were developed, including the variational multi-scale method (VMS) by Hughes {\it et al} \cite{Hughes1998}. In VMS, one casts simultaneously a coarse grid and a fine grid over the whole space domain, and describes a physical quantity in two scales, i.e., an average part identified in the coarse grid, and a fluctuation part in the fine grid. The fine scale fluctuation is reproduced from the coarse scale average by an operator based on Green's function. VMS was then used for large-eddy simulation of turbulent flows \cite{Codina2002}. An isogeometric analysis version was developed for more accurate geometrical representation \cite{Bazilevs2007}. Algebraic VMS formulations and adaptive strategies were developed to improve computational efficiency and accuracy \cite{Calo2014, Evans2013}. More recently, Leonor and Wagner extended VMS to solve heat conduction in LBPF, where fine-scale grid only covers a small subdomain around the heat source to resolve large temperature gradient \cite{Wagner2024}.  Making use of GPU acceleration and dynamic mode decomposition with control, the open-source GO-MELT framework provides a practical way for efficient thermal simulations and real-part qualification.

On the other hand, Tensor Decomposition (TD), or Proper Generalized Decomposition (PGD) when computed incrementally, is a model order reduction technique. It uses the sum of several modes, each in a variable separated form, to represent a function in much reduced degrees of freedom (DoFs). For elliptic partial differential equations, in particular the Poisson equation, Ammar {\it et al} \cite{Ammar2006} solved by PGD with alternating directions. Lu {\it et al} \cite{Lu2021} and Lei {\it et al} \cite{Zhang2022} proposed HiDeNN-PGD and HiDeNN-TD methods, respectively, by combining hierarchical deep-learning neural networks \cite{Zhang2021} to optimize both mesh geometry and the unknown function. Munoz {\it et al} \cite{Munoz2024} applied PGD with optimal transport for parametric analysis of the Poisson equation. For parabolic equations, in particular the heat equation, Ammar {\it et al} implemented PGD with a Galerkin variational formulation \cite{Ammar2007}. Bonithon {\it et al} combined PGD with the boundary element method, proposing a non-incremental solution strategy \cite{Bonithon2011}. Kpogan {\it et al} \cite{Kpogan2018} used PGD with the method of fundamental solutions. Tang {\it et al} \cite{Tang2023} introduced a preconditioner to accelerate PGD convergence. Recently Guo {\it et al} used TD to perform a large scale computation for heat transfer in additive manufacturing \cite{Guo2025}. For hyperbolic equations, Tang et al. proposed a way to expedite convergence for tensor decomposition in solving inviscid Burgers' equation and the Euler equations \cite{Tang2025,Xu2025}. For the advection-diffusion equation, which belongs to mixed hyperbolic-parabolic type, Nouy proposed a minimax PGD approach \cite{Nouy2010}, and Guan {\it et al} combined PGD with coordinate transformations \cite{Guan2024}. TD/PGD methods have been applied in viscoelasticity \cite{Ladeveze2010}, crack propagation \cite{Garikapati2020}, polymeric membrane \cite{Alotto2011}, microelectronics \cite{Krimm2019}, uncertainty quantification \cite{Wei2022}, nonlinear frequency response \cite{Lee2023}, and topology optimization \cite{Li2023,Hu2025}, etc.

To incorporate TD in VMS formulation, one may either take a space-time approach to decompose all space variables and the time variable, or a time-marching approach to decompose only the space variables while iterate in time. Though the former approach yields even smaller DoFs and avoids accumulative errors, numerically we observe difficulty in convergence. Notice that the total numerical cost is the product of number of iterations required for convergence (which is bigger than time-marching), and algebraic solver cost per iteration (which is smaller than time-marching). In contrast, the time-marching approach solves at each time step an elliptic equation with convergence guaranteed.

In this paper, we propose a multi-level VMS-TD algorithm for simulating heat conduction with moving heat source. In this algorithm, the temperature at both coarse-scale and fine-scale is described in TD (variable-separated) form, with fine scale part only resolved in a small subdomain that covers the heat source. Since the heat source moves, so does the subdomain. To tackle with the challenge in seamlessly combining the two scales together, first we take a second-order time semidiscretization. Secondly, at each time step, taking a two-level two space dimensional case as an example, we formulate the Galerkin weak form for the resulting coupled system with both fine and coarse scales. Thirdly, the coupled system is discretized and put in TD form. To attain the full competence of TD, it is critical to represent all involved components by arrays only, in particular, for the moving fine-scale subdomain that covers the heat source, and the interface across scales. So fourthly, we design special treatments to data in the moving subdomain and the interface. With all these, we are ready to construct a two-level two space dimensional VMS-TD algorithm. It can be readily extended to multi-level and high space dimensional simulations. Compared to existing simulation frameworks including GO-MELT, the TD strategy substantially reduces DoFs and consequently the computational time. Moreover, VMS-TD has second order accuracy in time, as compared to GO-MELT which uses first order forward Euler scheme for time integration. We also note that this work performs TD only in space dimensions, whereas a fully space-time decomposed TD has been developed in another recent work allowing higher and tunable order of accuracy  \cite{Zhanglei2025}.

The rest of this paper is organized as follows. In Section \ref{sect: basic info}, we formulate the multi-level VMS-TD algorithm for solving the heat equation.
In Section \ref{sect: Seamless treatment combining VMS with TD}, we describe in detail the special treatments for interface and moving subdomain, and present the algorithm flowcharts. In Section \ref{sect: numerical results}, we demonstrate the effectiveness of the VMS-TD algorithm through a two-level two-dimensional example, and a three-level three-dimensional example. Then we compare computational costs by contrasting the proposed algorithm with a fully resolved fine scale TD simulation for three-level three-dimensional heat conduction with fixed heat source. Finally, in Section \ref{sect: conclusion}, we make some concluding remarks. A preliminary multiple time step algorithm is also given in Appendix.

\section{VMS-TD formulation for heat equation} \label{sect: basic info}

In this section, we propose a TD algorithm to implement VMS formulation to the heat equation.

\subsection{Heat equation and time-discretization} \label{sect: heat equ.}

Consider heat conduction in $\Omega \subset \mathbb{R}^d$ for time from $0$ to $T$, namely
\begin{equation} \label{eq: heat equ.}
	u_t - \nu \Delta u = f, \quad x \in \Omega, t \in [0, T].
\end{equation}
Here $u(x, t)$ is the temperature, $\nu > 0$ the diffusivity, $\Delta = \partial_{x_1}^2 + \cdots + \partial_{x_d}^2$ the Laplace operator, $f$ a moving heat source with compact support. Without loss of generality, let $\Omega = [0, L_1]^d$ ($d=2, 3$). Initial and boundary conditions are
\begin{equation}
	u(x, 0) = u_0(x), \quad \left. u(x, t) \right|_{\partial \Omega} = u_B,
\end{equation}
with the ambient temperature $u_B = 0$.

Taking a time step size $\Delta t = T / N_t$, we denote $t_n = n \cdot \Delta t$, and $u^n(x) = u(x, t_n)$. A second order difference semi-discretization reads
\begin{equation} \label{eq: time-discretized heat equ.}
	\dfrac{u^n - u^{n-1}}{\Delta t} - \nu \Delta \dfrac{u^n + u^{n-1}}{2} = f^{n - 1/2}, \quad n = 1, \cdots, N_t,
\end{equation}
with $f^{n-1/2}(x) = f(x, (t_n + t_{n-1})/2)$. The Galerkin weak form reads
\begin{equation} \label{eq: time-discretized heat equ. G. weak form}
	a(w^n, u^n)_{\Omega} = \hat{a}(w^{n}, u^{n-1})_{\Omega} + (w^n, f^{n-1/2})_{\Omega}, \quad \forall w^n \in H_0^1(\Omega),
\end{equation}
where
\begin{equation} \label{eq: def operator a}
	a(w, u)_\Omega = \dint_{\Omega} \left( \dfrac{1}{\Delta t} w  u + \dfrac{\nu}{2} \nabla w \cdot \nabla u \right) \ud x ,
\end{equation}
\begin{equation} \label{eq: def operator a_hat}
	\hat{a}(w, u)_\Omega = \dint_{\Omega} \left( \dfrac{1}{\Delta t} w  u - \dfrac{\nu}{2} \nabla w \cdot \nabla u \right) \ud x ,
\end{equation}
\begin{equation} \label{eq: def operator inner product of function}
	(w, f)_\Omega = \dint_{\Omega} w  f \ud x ,
\end{equation}
and $\nabla = [\partial_{x_1}, \cdots, \partial_{x_d}]$ is the gradient operator, $\cdot$ denotes inner product for vectors.

% \section{Two-level model} \label{sect: two-level}

\subsection{VMS formulation} \label{sect: 2-lv:VMS}

Following \cite{Hughes1998, Wagner2024}, we assign a subdomain $\Theta^n$ for each time step $t^n$. This subdomain covers the heat source at this step, therefore moves along with the heat source. $u^n(x, y)$ is decomposed into its mean part
\begin{equation}
	\overline{u^n}(x, y) = \mathcal{E}(u^n(x, y) ),
\end{equation}
with $\mathcal{E}$ denoting an averaging operator, and fluctuation part
\begin{equation} \label{eq: interscale constraint}
	\widetilde{u^n}(x, y) = u^n(x, y) - \overline{u^n}(x, y).
\end{equation}
Since there is no fine grid out of $\Theta^n$ at $t^n$ step, we take
\begin{equation} \label{eq: u' = 0 outside Theta}
	\widetilde{u^n}(x, y) = 0, \quad (x, y) \in \Omega \backslash \Theta^n.
\end{equation}

In $\Theta^n$, we consider the fluctuation part of an arbitrary test function $\widetilde{w^n}$, and obtain the weak form
\begin{equation} \label{eq: VMS fine}
	a(\widetilde{w^n}, u^n)_{\Theta^n} = \hat{a}(\widetilde{w^n}, u^{n-1})_{\Theta^n} + (\widetilde{w^n}, f^{n-1/2})_{\Theta^n},
\end{equation}
with artificial boundary condition
\begin{equation} \label{eq: VMS.BC fine}
	 u^n(x, y)  = \overline{u^n}(x, y), \quad (x, y) \in \partial\Theta^n.
\end{equation}

Over the coarse grid, we consider the mean part of the test function $\overline{w^n}$. The weak form reads
\begin{equation} \label{eq: VMS coarse}
	a(\overline{w^n}, u^n)_{\Omega} = \hat{a}(\overline{w^n}, u^{n-1})_{\Omega} + (\overline{w^n}, f^{n-1/2})_{\Omega},
\end{equation}
with physical boundary condition
\begin{equation} \label{eq: VMS.BC coarse}
	\overline{u^n}(x, y) = u_{B}, \quad (x, y) \in {\partial\Omega}.
\end{equation}
Due to the assumption \eqref{eq: u' = 0 outside Theta}, we rewrite
\begin{equation}
	a(\overline{w^n}, u^n)_{\Omega} = a(\overline{w^n}, \overline{u^n})_{\Omega} + a(\overline{w^n}, u^n )_{\Theta^n} - a(\overline{w^n}, \mathcal{E}(u^n))_{\Theta^n},
\end{equation}
and similar for $a(\overline{w^n}, u^{n-1})_\Omega$. The weak form amounts to
\begin{equation} \label{eq: VMS coarse, eliminate u'}
	\begin{array}{rl}
	& a(\overline{w^n}, \overline{u^n})_{\Omega} + a(\overline{w^n}, u^n )_{\Theta^n} - a(\overline{w^n}, \mathcal{E}(u^n))_{\Theta^n} \\[3mm]
	= & \hat{a}(\overline{w^n}, \overline{u^{n-1}})_{\Omega}
	+ \hat{a}(\overline{w^n}, u^{n-1})_{\Theta^{n-1}}
	- \hat{a}(\overline{w^n}, \mathcal{E}(u^{n-1}))_{\Theta^{n-1}}
	+ (\overline{w^n}, f^{n-1/2})_{\Omega}.
	\end{array}
\end{equation}
Equations \eqref{eq: VMS fine}\eqref{eq: VMS.BC fine}\eqref{eq: VMS.BC coarse}\eqref{eq: VMS coarse, eliminate u'} form a coupled system for unknown variable $\overline{u^n}(x, y)$ in $\Omega$ and $u^n(x, y)$ in $\Theta^n$.

\subsection{Finite element analysis formulation}

\begin{figure}[H]
	\centering
	\includegraphics[width=.5\textwidth]{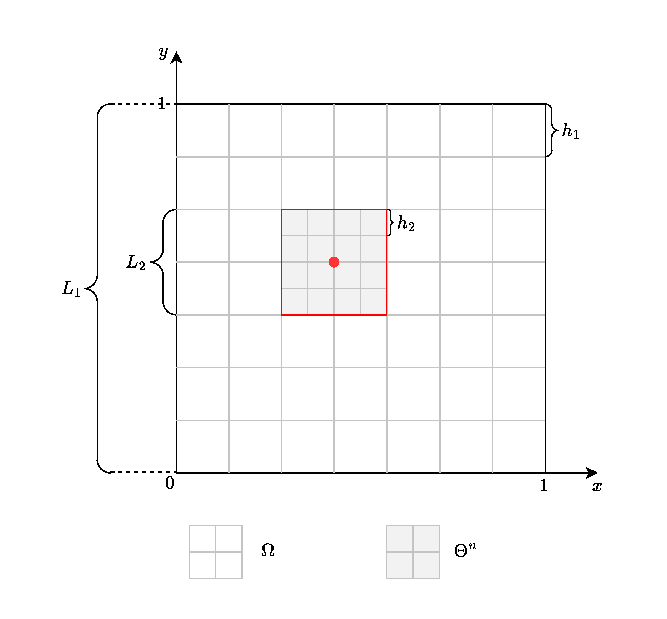}
	\caption{Schematical plot of a VMS in 2D.}
	\label{fig: VMS 2D 2L mesh}
\end{figure}

Cast a uniform coarse mesh with mesh size $h_1 = \dfrac{L_1}{N_1}$ in each dimension of $\Omega$, and refine it with mesh size $h_2 = \dfrac{L_2}{N_2}$ in each dimension of $\Theta^n$. Fig. \ref{fig: VMS 2D 2L mesh} shows a 2D mesh.

Using vector index $\bm{J} = (J_x, J_y)$ for coarse mesh nodes, grid function $\overline{U^n_{\bm{J}}} \sim \overline{u^n}(x_{J_x}, y_{J_y})$, and shape function $N^1_{\bm{J}}(x, y)$, we take finite element interpolation
\begin{equation} \label{eq: interpolation coarse}
	\overline{u^n}(x, y) = \dsum_{\bm{J}} N^1_{\bm{J}}(x, y) \overline{U^n_{\bm{J}}}  , \quad (x, y) \in \Omega.
\end{equation}

Similarly, using vector index $\bm{j} = (j_x, j_y)$ for fine mesh nodes in $\Theta^n$, grid function $U^n_{\bm{j}} = u^n(x_{j_x}, y_{j_y})$, and shape function $N_{\bm{j}}^2(x, y)$, we take finite element interpolation
\begin{equation} \label{eq: interpolation fine}
	u^n(x, y) = \dsum_{\bm{j}}  N^2_{\bm{j}}(x, y) U^n_{\bm{j}}, \quad (x, y) \in \Theta^n.
\end{equation}

By standard finite element analysis, we obtain the following discrete system.
\begin{subnumcases}{\label{eq: discrete VMS system}}
	\begin{array}{l}
		\begin{array}{l}
			\quad \displaystyle  L^{11.\Omega} \overline{U^n} +  L^{12.\Theta^n} U^n - L^{11.\Theta^n} [\mathcal{E}(U^n)] \\[3mm]
			 \displaystyle =  \hat{L}^{11.\Omega} \overline{U^{n-1}} + \hat{L}^{12.\Theta^{n-1}} U^{n-1} -  \hat{L}^{11.\Theta^{n-1}} [\mathcal{E}(U^{n-1})]
			+ F^{1.\Omega};%, \quad \text{Inner points};
			\end{array}
	\end{array} \label{eq: discrete VMS coarse: equ.} \\
	\begin{array}{l}
			\displaystyle L^{22.\Theta^n} U^n
			\displaystyle =  \hat{L}^{22.\Theta^n} U^{n-1}
			+ F^{2.\Theta^n};
		\end{array} \label{eq: discrete VMS fine: equ.} \\
		\overline{U^n} = U_B^n; \label{eq: discrete VMS coarse: b.c.} \\
		U^n = \mathcal{I}(\overline{U^n}). \label{eq: discrete VMS fine: b.c.}
\end{subnumcases}
Here the matrix superscript identifies the scales and domains. For instance, the matrix $L^{12.\Theta^n}$ arises from the bilinear form $a(\overline{w^n}, u^n)_{\Theta^n}$ in \eqref{eq: VMS coarse, eliminate u'}, with $12$ identifying coarse scale (scale 1) variable $\overline{w^n}$ and fine scale (scale 2) variable $u^n$, and $\Theta^n$ identifying the domain for integration. Its $(\bm{J}, \bm{j})$ entry is given by
\begin{equation} \label{eq: def L operator}
		\displaystyle L^{12.\Theta^n}_{\bm{J}\bm{j}} = \dfrac{1}{\Delta t} \int_{\Theta^n}  N_{\bm{J}}^1(x, y) N_{\bm{j}}^2(x, y) \ud (x, y)
		\displaystyle + \dfrac{\nu}{2} \int_{\Theta^n}  \nabla N_{\bm{J}}^1(x, y) \cdot \nabla N_{\bm{j}}^2(x, y)  \ud (x, y).
\end{equation}
All other terms are similarly defined.

The averaging operator $\mathcal{E}$ maintains the same meaning, though it applies to the vector $U^n$. An interpolation operator $\mathcal{I}$ is introduced on $\partial \Theta^n$ for interpolating the coarse grid nodal values to provide fine scale ones.

\subsection{Tensor decomposition formulation} \label{sect: tensor decomposition formulation}
In tensor decomposition method, one represents an unknown by summing a number of modes in variable-separated form. Taking the fine scale temperature $U^n$ (an $(N_2+1) \times (N_2+1)$ matrix) as an example, we approximate by column vectors $X_1^{2,n}, Y_1^{2,n}, \cdots, X_{Q_2}^{2,n}, Y_{Q_2}^{2,n}$, each of order $N_2$. Here '2' again indicates the fine scale.
\begin{equation}
	U^n_{\bm{j}} \approx \sum_{q=1}^{Q_2} X_{q,j_x}^{2,n} Y_{q,j_y}^{2,n}, \quad j_x, j_y = 0, \cdots N_2.
\end{equation}
Similarly, we approximate the coarse scale temperature $\overline{U^n}$ by
\begin{equation}
	\overline{U^n_{\bm{J}}} \approx \sum_{q=1}^{Q_1} X_{q,J_x}^{1,n} Y_{q,J_y}^{1,n}, \quad J_x, J_y = 0, \cdots N_1.
\end{equation}
Here '1' again indicates the fine scale. It is natural to take variable separated shape functions
\begin{equation}
	N_{\bm{J}}^1(x, y) = N^1_{J_x}(x) N^1_{J_y}(y), \quad N_{\bm{j}}^2(x, y) = N^2_{j_x}(x) N^2_{j_y}(y).
\end{equation}
Then the matrices in \eqref{eq: discrete VMS system} take a similar form, e.g.,
\begin{equation} \label{eq: TD formed L operator}
	\begin{array}{rl}
		\displaystyle L^{12.\Theta^n}_{\bm{J}\bm{j}} = & \dfrac{1}{\Delta t} \dint N_{J_x}^1(x)N_{j_x}^2(x)  \ud x \dint N_{J_y}^1(y)  N_{j_y}^2(y)  \ud y  \\[3mm]
		& + \dfrac{\nu}{2} \int \dfrac{\ud}{\ud x} N_{J_x}^1(x) \dfrac{\ud}{\ud x} N_{j_x}^2(x) \ud x \dint N_{J_y}^1(y) N_{j_y}^2(y) \ud y   \\[3mm]
		& + \dfrac{\nu}{2} \int N_{J_x}^1(x) N_{j_x}^2(x) \ud x \dint \dfrac{\ud}{\ud y} N_{J_y}^1(y) \dfrac{\ud}{\ud y} N_{j_y}^2(y) \ud y   \\[3mm]
		\equiv & \dsum_{\alpha=1}^{3} L^{12.\Theta^n}_{x, \alpha J_x j_x} L^{12.\Theta^n}_{y, \alpha J_y j_y}.
	\end{array}
\end{equation}
Here the integrations in $x$ and $y$ are taken over the interval in each dimension of $\Theta^n$.

All other matrices are expressed and abbreviated in a similar fashion. The source term usually needs more modes to approximate, e.g.,
\begin{equation}
	F^{1.\Omega}_{\bm{J}} = \sum_{\beta=1}^{Q_f} F^{1.\Omega}_{x, \beta J_x} F^{1.\Omega}_{y, \beta J_y}, \quad J_x, J_y = 0, \cdots, N_1,
\end{equation}
\begin{equation}
	F^{2.\Theta^n}_{\bm{j}} = \sum_{\beta=1}^{Q_f} F^{2.\Theta^n}_{x, \beta j_x} F^{2.\Theta^n}_{y, \beta j_y}, \quad j_x, j_y = 0, \cdots, N_2.
\end{equation}
In our numerical simulations, we always take boundary condition to be zero.
Non-zero boundary conditions may be treated as source terms.

Each term in the finite element analysis formulation \eqref{eq: discrete VMS system} may now be expressed in a dimension separated form, e.g., the $\bm{J}$-th entry
\begin{equation}
	\left( L^{12.\Theta^n} U^n \right)_{\bm{J}} = \dsum_{\alpha = 1}^{3} \sum_{q = 1}^{Q_2}
	\left( \sum_{j_x = 0}^{N_2} L^{12.\Theta^n}_{x, \alpha J_x j_x} X^{2,n}_{q, j_x} \right)
	\left( \sum_{j_y = 0}^{N_2} L^{12.\Theta^n}_{y, \alpha J_y j_y} X^{2,n}_{q, j_y} \right).
\end{equation}

Except for the treatment of averaging $\mathcal{E}$, interpolation $\mathcal{I}$ and possibly $\Theta^n \neq \Theta^{n-1}$, the finite element analysis formulation \eqref{eq: discrete VMS system} may be formally written as
\begin{subnumcases}{\label{eq: discrete VMS-TD system}}
	\dsum_{\alpha = 1}^{A_{11}} \dsum_{q=1}^{Q_1} \bm{L}^{11}_{x, \alpha} \bm{X}^{1,n}_{q} \otimes \bm{L}^{11}_{y, \alpha} \bm{Y}^{1,n}_{q}
	+ \dsum_{\alpha = 1}^{A_{12}} \dsum_{q=1}^{Q_2} \bm{L}^{12}_{x, \alpha} \bm{X}^{2,n}_{q} \otimes \bm{L}^{12}_{y, \alpha} \bm{Y}^{2,n}_{q}
	= \dsum_{\beta = 1}^{B_1} \bm{F}^1_{x, \beta} \otimes \bm{F}^1_{y, \beta} ; \label{eq: discrete VMS-TD system: coarse} \\
	\dsum_{\alpha = 0}^{A_{21}} \dsum_{q=1}^{Q_1} \bm{L}^{21}_{x, \alpha} \bm{X}^{1,n}_{q} \otimes \bm{L}^{21}_{y, \alpha} \bm{Y}^{1,n}_{q}
	+ \dsum_{\alpha = 0}^{A_{22}} \dsum_{q=1}^{Q_2} \bm{L}^{22}_{x, \alpha} \bm{X}^{2,n}_{q} \otimes \bm{L}^{22}_{y, \alpha} \bm{Y}^{2,n}_{q}
	= \dsum_{\beta = 2}^{B_2} \bm{F}^2_{x, \beta} \otimes \bm{F}^2_{y, \beta} ; \label{eq: discrete VMS-TD system: fine}
\end{subnumcases}
Here $ \bm{X}^{1,n}_{q}$ is the matrix formed with $X^{1,n}_{q, 0}, \cdots, X^{1,n}_{q, N_x^1}$.
Similarly define other matrices in \eqref{eq: discrete VMS-TD system}.
We adopt an iterative scheme to solve this system. Take the first set of equations as an example,
\begin{equation}
	\dsum_{\alpha = 1}^{A_{11}} \dsum_{q=1}^{Q_1} \bm{L}^{11}_{x, \alpha} \bm{X}^{1,n}_{q} \otimes \bm{L}^{11}_{y, \alpha} \bm{Y}^{1,n}_{q}
	= \dsum_{\beta = 1}^{B_1} \bm{F}^1_{x, \beta} \otimes \bm{F}^1_{y, \beta}
	- \dsum_{\alpha = 1}^{A_{12}} \dsum_{q=1}^{Q_2} \bm{L}^{12}_{x, \alpha} \bm{X}^{2,n}_{q} \otimes \bm{L}^{12}_{y, \alpha} \bm{Y}^{2,n}_{q}
	\equiv \dsum_{\beta = 1}^{\hat{B}_1} \hat{\bm{F}}^1_{x, \beta} \otimes \hat{\bm{F}}^1_{y, \beta},
\end{equation}
or component-wise,
\begin{equation} \label{eq: general discrete equation: simplified form}
	\dsum_{\alpha = 1}^{A_{11}} \dsum_{q=1}^{Q_1}
	\left( \dsum_{J_x = 0}^{N_1} {L}^{11}_{x, \alpha R_x J_x} {X}^{1,n}_{q, J_x} \right)
	\left( \dsum_{J_y = 0}^{N_1}  {L}^{11}_{y, \alpha R_y J_y} {Y}^{1,n}_{q, J_y} \right)
	= \dsum_{\beta = 1}^{\hat{B}_1} \hat{{F}}^1_{x, \beta R_x} \hat{{F}}^1_{y, \beta R_y}, \quad R_x, R_y = 0, \cdots, N_1.
\end{equation}

First, multiply \eqref{eq: general discrete equation: simplified form} by $Y^{1,n}_{i,R_y}$ and sum over $R_y$ to get a linear system about $X^{1,n}_{q,J_x}$,
\begin{equation} \label{eq: general discrete equation a.d.s. x dir.}
	\begin{array}{l}
		\quad \displaystyle \sum_{q=1}^{Q_1} \dsum_{J_x = 0}^{N_1} \left[ \sum_{\alpha=1}^{A_{11}} \left( \dsum_{R_y = 0}^{N_1} \dsum_{J_y = 0}^{N_1}  Y^{1,n}_{i,R_y} {L}^{11}_{y, \alpha R_y J_y} Y^{1,n}_{q, J_y} \right) {L}^{11}_{x, \alpha R_x J_x}  \right] X^{1,n}_{q,J_x}  \\
		\displaystyle = \sum_{\beta=1}^{\hat{B}_1} \hat{{F}}^1_{x, \beta R_x} \left( \dsum_{R_y = 0}^{N_1} \hat{{F}}^1_{y, \beta R_y} Y^{1,n}_{i,{R_y}} \right)
		\displaystyle, \quad i = 1, \cdots, Q_1, \quad R_x = 0, \cdots, N_1.
		\end{array}
\end{equation}
Then multiply \eqref{eq: general discrete equation: simplified form} by $X^{1,n}_{i,R_x}$ and sum over $R_x$ to get a linear system about $Y^{1,n}_{q,J_y}$,
\begin{equation} \label{eq: general discrete equation a.d.s. y dir.}
	\begin{array}{l}
		\quad \displaystyle \sum_{q=1}^{Q_1} \sum_{J_y=0}^{N_1} \left[ \sum_{\alpha=1}^{A_{11}} \left( \dsum_{R_x = 0}^{N_1} \dsum_{J_x = 0}^{N_1} X^{1,n}_{i,R_x} {L}^{11}_{x, \alpha R_x J_x} X^{1,n}_{q,J_x} \right)  {L}^{11}_{y, \alpha R_y J_y}  \right]  Y^{1,n}_{q,J_y} \\
		\displaystyle
		= \sum_{\beta=1}^{\hat{B}_1} \left( \sum_{R_x=0}^{N_1} \hat{{F}}^1_{x, \beta R_x} X^{1,n}_{i,R_x} \right) \hat{{F}}^1_{y, \beta R_y}
		\displaystyle, \quad i = 1, \cdots, Q_1, \quad R_y = 0, \cdots, N_1.
		\end{array}
\end{equation}
To ensure the uniqueness, we add a constraint
\begin{equation} \label{eq: X_norm = 1}
	\sum_{q=1}^{Q_1} \sum_{J_x = 0}^{N_1} \left(X^{1,n}_{q,J_x}\right)^2 = 1.
\end{equation}
Iterate alternately \eqref{eq: general discrete equation a.d.s. x dir.} and \eqref{eq: general discrete equation a.d.s. y dir.} till convergence, then we get the solution of \eqref{eq: discrete VMS-TD system: coarse}. Similarly we solve \eqref{eq: discrete VMS-TD system: fine}. In a loop we iterate alternately \eqref{eq: discrete VMS-TD system: coarse} and \eqref{eq: discrete VMS-TD system: fine} till convergence, then the coupled system \eqref{eq: discrete VMS-TD system} is solved.

\section{VMS-TD Algorithm} \label{sect: Seamless treatment combining VMS with TD}

The VMS method significantly reduces DoFs to $(N_1)^d + (N_2)^d$. Nevertheless, for large $N_1$ and $N_2$, this remains demanding. By using the TD algorithm, we further reduce DoFs to $(N_1 Q_1 + N_2 Q_2) d$, where $Q_1$ and $Q_2$ are the number of TD modes on the coarse and fine scales, respectively.

To attain the numerical efficiency offered by VMS and TD, it is crucial to avoid reconstructing data.
This requires to deal with the interface between the two scales, and the moving fine-scale subdomain.
In the following, we propose special treatments to ensure that all data are stored in arrays.

\subsection{Interface} \label{sect: interface boundary}

\begin{figure}[H]
    \centering
    \includegraphics[width=\textwidth]{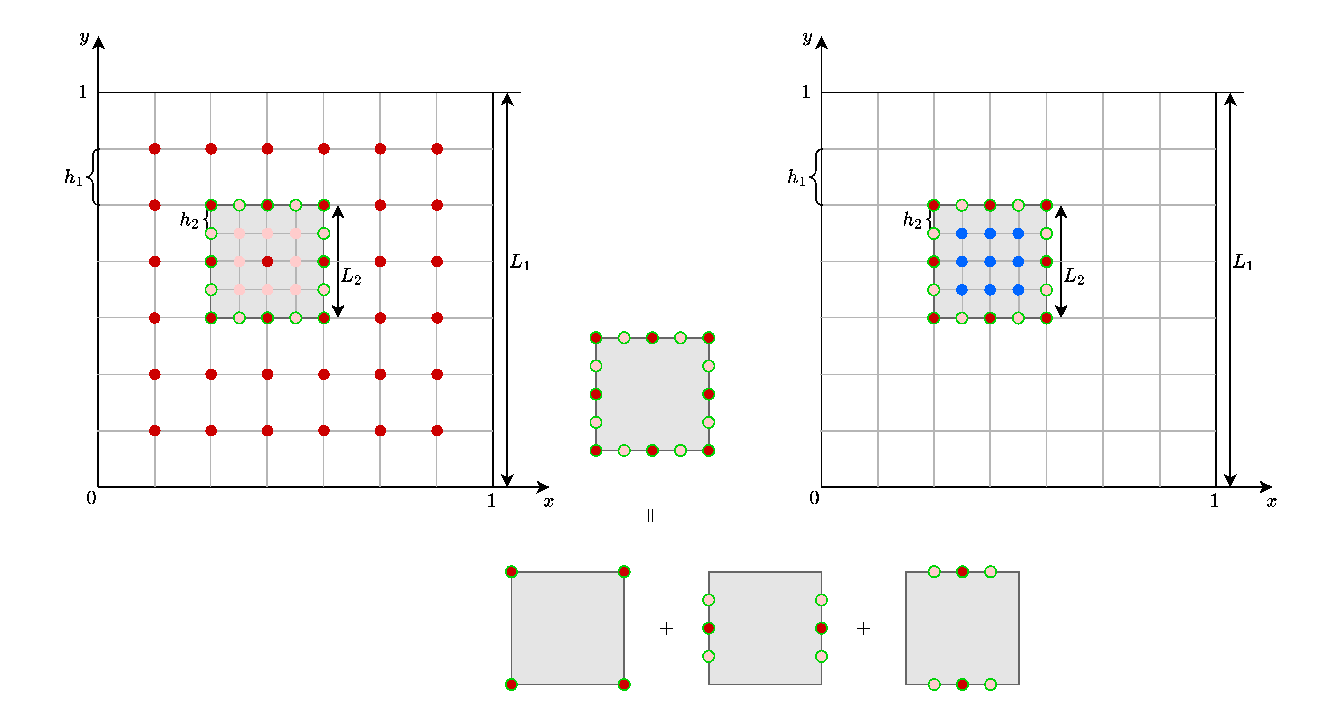}
    \caption{Interface treatment. upper-left: coarse grid on the whole domain; upper-right: fine grid; upper-middle: the interface; lower-middle: three types of grid points at the interface.}
    \label{fig: treatment for interface boundary}
\end{figure}

We assign the boundary condition on $\partial \Theta^n$ by interpolating the coarse scale temperature $\overline{U^n}$. It is desirable to form the boundary condition as the sum of variable-separated modes. We schematically illustrate this in Fig. \ref{fig: treatment for interface boundary}.
Grid points marked by green circle represent the interface across two scales.
They are on the fine-scale subdomain.
This shape does not support variable-separated functions.
We divide such points into three types: corner points, vertical boundary, and horizontal boundary, as shown at the lower-subplot of Fig. \ref{fig: treatment for interface boundary}.
Each of them supports variable-separated functions with one-dimensional data.

% The mathematical formulation of the interface treatment is as follows.
Let the coarse scale grid points in $\Theta^n$ be indexed by $J_x = J_x^{\text{min}}, \cdots, J_x^{\text{max}}$; and $J_y = J_y^{\text{min}}, \cdots, J_y^{\text{max}}$. Corresponding coarse grid nodal values are
\begin{equation}
	\overline{U^n}_{J_x J_y} = \sum_{q = 1}^{Q_1} X^{1,n}_{q,J_x} Y^{1,n}_{q,J_y}.
\end{equation}

The shape function $N^1_{J_x}(x) N^1_{J_y}(y)$ yields interpolated temperature at fine scale nodes
\begin{equation}
	U^n_{j_x j_y} = \dsum_{q=1}^{Q_1} \dsum_{J_x = J_x^{\text{min}}}^{J_x^{\text{max}}} N^1_{J_x}(x^2_{j_x}) X^{1,n}_{q,J_x} \dsum_{J_y = J_y^{\text{min}}}^{J_y^{\text{max}}} N^1_{J_y}(y^2_{j_y}) Y^{1,n}_{q,J_y}
	\equiv \dsum_{q=1}^{Q_1} \left(\mathcal{I} X^{1,n}_{q}\right)_{j_x} \left(\mathcal{I} Y^{1,n}_{q}\right)_{j_y}.
\end{equation}
Note that the latter two summations on $J_x$, $J_y$ may be computed separately.

Since we are only interested in the interfacial nodes, we may consider $ \left. U^n \right|_{\partial \Theta^n}$ as the summation of three parts, each in variable-separated form.

The first part takes care of the vertices only, namely, expressing
\begin{equation} \label{eq: InterBd.MatForm.BB}
	\left[
	\begin{array}{ccccc}
	U^n_{00} & 0 & \cdots & 0 & U^n_{0 N_2} \\
	0 & 0 & \cdots & 0 & 0 \\
	\cdots & \cdots & \cdots & \cdots & \cdots \\
	0 & 0 & \cdots & 0 & 0 \\
	U^n_{N_2 0} & 0 & \cdots & 0 & U^n_{N_2 N_2} \\
	\end{array}
	\right]
	= \dsum_{q=1}^{Q_1}
	\left[
	\begin{array}{c}
	\left(\mathcal{I} X^{1,n}_{q}\right)_{0} \\ 0 \\ \cdots \\ 0 \\ \left(\mathcal{I} X^{1,n}_{q}\right)_{N_2}
	\end{array}
	\right]
	\left[
	\begin{array}{ccccc}
	\left(\mathcal{I} Y^{1,n}_{q}\right)_{0} & 0 & \cdots & 0 & \left(\mathcal{I} Y^{1,n}_{q}\right)_{N_2} \\
	\end{array}
	\right].
\end{equation}

The second part takes care of the two vertical edges (except the vertices),
\begin{equation}\label{eq: InterBd.MatForm.BI}
	\left[
	\begin{array}{ccccc}
	0 & U^n_{01} & \cdots & U^n_{0, N_2-1} & 0 \\
	0 & 0 & \cdots & 0 & 0 \\
	\cdots & \cdots & \cdots & \cdots & \cdots \\
	0 & 0 & \cdots & 0 & 0 \\
	0 & U^n_{N_2 1} & \cdots & U^n_{N_2, N_2-1} & 0 \\
	\end{array}
	\right]
	= \dsum_{q=1}^{Q_1}
	\left[
	\begin{array}{c}
	\left(\mathcal{I} X^{1,n}_{q}\right)_{0} \\ 0 \\ \cdots \\ 0 \\ \left(\mathcal{I} X^{1,n}_{q}\right)_{N_2}
	\end{array}
	\right]
	\left[
	\begin{array}{ccccc}
	0 & \left(\mathcal{I} Y^{1,n}_{q}\right)_{1} & \cdots & \left(\mathcal{I} Y^{1,n}_{q}\right)_{N_2-1} & 0 \\
	\end{array}
	\right].
\end{equation}

Similarly, the third part takes care of the two horizontal edges
\begin{equation}\label{eq: InterBd.MatForm.IB}
	\left[
	\begin{array}{ccccc}
	0 & 0 & \cdots & 0 & 0 \\
	U^n_{10} & 0 & \cdots & 0 & U^n_{N_2-1, 0} \\
	\cdots & \cdots & \cdots & \cdots & \cdots \\
	U^n_{1 N_2} & 0 & \cdots & 0 & U^n_{N_2-1, N_2} \\
	0 & 0 & \cdots & 0 & 0 \\
	\end{array}
	\right]
	= \dsum_{q=1}^{Q_1}
	\left[
	\begin{array}{c}
	0 \\ \left(\mathcal{I} X^{1,n}_{q}\right)_{1} \\ \cdots \\ \left(\mathcal{I} X^{1,n}_{q}\right)_{N_2-1} \\ 0
	\end{array}
	\right]
	\left[
		\begin{array}{ccccc}
		\left(\mathcal{I} Y^{1,n}_{q}\right)_{0} & 0 & \cdots & 0 & \left(\mathcal{I} Y^{1,n}_{q}\right)_{N_2} \\
		\end{array}
	\right].
\end{equation}
Summing them together, we get a tensor decomposed form of the boundary interpolation. In implementation, we may further suppress the number of modes ($Q_1$) by singular value decomposition.

We store data in arrays for the right hand side of \eqref{eq: InterBd.MatForm.BB}\eqref{eq: InterBd.MatForm.BI}\eqref{eq: InterBd.MatForm.IB}. For a general $d$-dimensional case, there are $(2^d - 1)Q_1$ modes. Although the number of modes increases, all data are stored in arrays.

\subsection{Moving fine scale subdomain} \label{sect: moving region}

\begin{figure}[H]
    \centering
    \includegraphics[width=.5\textwidth]{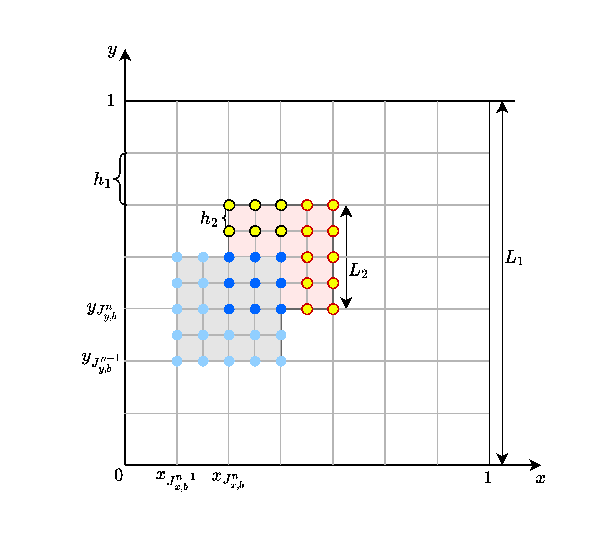}
    \caption{Moving fine-scale subdomain. blue: $U^{n-1}$; dark-blue: reserved data; light-blue: discarded data; yellow: re-built data in $\Theta^n$; red- and black-circled: occupy rectangular regions allowing variable-separated form.}
    \label{fig: treatment for moving region}
\end{figure}

As the heat source moves, the subdomain moves accordingly.
Using the known data $U^{n-1}$ in $\Theta^{n-1}$, we re-build data in $\Theta^n$ to supply the governing equation \eqref{eq: discrete VMS fine: equ.}.
See Fig. \ref{fig: treatment for moving region} for illustration, where the subdomain moves to upper-right.
Denoting the lower-left vertex of $\Theta^n$ as $\left( x_{ J_{x,b}^n }, y_{ J_{y,b}^n } \right)$, and that of $\Theta^{n-1}$ as $\left( x_{ J_{x,b}^{n-1} }, y_{ J_{y,b}^{n-1} } \right)$, we shall discard $U^{n-1}$ in the '$L$'-shaped subdomain $\Theta^{n-1} \backslash \Theta^{n}$ (light-blue points), and supplement $U^{n-1}$ in the '$7$'-shaped subdomain $\Theta^{n} \backslash \Theta^{n-1}$ (yellow points).

% It is mathematically clean if we do not consider tensor decomposition.
% It is not natural for variable-separated form to express a function with a '$7$'-shaped support when programming.
To this end, we split yellow points into two types. One type is represented by red-circled points, and the other by black-circled points.
Both of them occupy rectangular regions, allowing variable-separated form of each mode.

% When the fine mesh region moves (i.e., $\Theta^n \neq \Theta^{n-1}$), data in $U^{n-1}$ in \eqref{eq: discrete VMS fine: equ.} must be shifted in the opposite direction, as shown in Figure \ref{fig: treatment for moving region}.

\subsection{Numerical algorithm}

Now we are ready to describe the overall algorithm. For the sake of clarity, we omit some of the details such as the above treatments for interfaces and the moving fine scale subdomain.

For one-step marching from $t_{n-1}$ to $t_n$ of \eqref{eq: VMS fine}\eqref{eq: VMS.BC fine}\eqref{eq: VMS.BC coarse}\eqref{eq: VMS coarse, eliminate u'}, we have $U^{n-1}$ on $\Theta^{n-1}$ and $\overline{U^{n-1}}$ on $\Omega$. At step $t_n$, the fine-scale subdomain
\begin{equation}
	\Theta^n = \left\{ (x, y) \left|  x_c^n - \dfrac{L_2}{2} < x < x_c^n + \dfrac{L_2}{2}, \quad y_c^n - \dfrac{L_2}{2} < y < y_c^n + \dfrac{L_2}{2} \right. \right\},
\end{equation}
where the center $\left( x_c^n, y_c^n \right)$ is the node closest to the moving heat source. With a slight abuse of notation, we still take $U^{n-1}$ to denote the value on $\Theta^n$ which actually has been treated in afore-mentioned manner.

We take an iterative scheme to solve \eqref{eq: discrete VMS system}, and an alternating direction scheme for TD in each iteration. More precisely, as described in subsection \ref{sect: tensor decomposition formulation}, we initialize $X^{1,n}_{q, J_x} = \sin (q \pi J_x / N)$, $Y^{1,n}_{q, J_y} = \sin (q \pi J_y / N)$, $U^n = 0$, and TD algorithm to solve $\overline{U^n}$ from \eqref{eq: discrete VMS coarse: equ.}\eqref{eq: discrete VMS coarse: b.c.}. Next, using this $\overline{U^n}$, we solve $U^n$ from \eqref{eq: discrete VMS fine: equ.}\eqref{eq: discrete VMS fine: b.c.} by TD algorithm. In our implementation, the convergence criterion between the coarse and fine scales is set as
\begin{equation} \label{eq: criteria for inter-scale iteration}
    \dfrac{ \left\| \overline{U^{n,\theta}} - \overline{U^{n,\theta - 1}} \right\|}{ \left\| \overline{U^{n,\theta}} \right\|} < 10^{-3},
\end{equation}
with $\theta$ the iteration count. The TD convergence criterion is set as
\begin{equation} \label{eq: criteria for tensor decomposition}
    \sqrt{ \dfrac{\left\Vert X^{1,n,\theta,\rho} - X^{1,n,\theta,\rho-1} \right\Vert^2}{\left\Vert X^{1,n,\theta,\rho} \right\Vert^2} + \dfrac{\left\Vert Y^{1,n,\theta,\rho} - Y^{1,n,\theta,\rho-1} \right\Vert^2}{\left\Vert Y^{1,n,\theta,\rho} \right\Vert^2} } < 10^{-2},
\end{equation}
with $\rho$ the iteration count. Note that the Frobenius norms are adopted for these criteria.

The above algorithm is summarized in a flowchart shown below.

\begin{algorithm}[H]
	\SetAlgoLined
	
	Initialize $U^{n, 0} \leftarrow 0$\;
	\For{$\theta = 1, \cdots, \theta_{\text{max}}$}{
		Formulate \eqref{eq: discrete VMS coarse: equ.} and \eqref{eq: discrete VMS coarse: b.c.} into the form of \eqref{eq: general discrete equation a.d.s. x dir.} and \eqref{eq: general discrete equation a.d.s. y dir.}\;
		Initialize $X^{1,n, \theta, 0}_{q, J_x} \leftarrow \sin (q \pi J_x / N)$, $Y^{1,n, \theta, 0}_{q, J_y} \leftarrow \sin (q \pi J_y / N)$\;
		\For{$\rho = 1, \cdots, \rho_{\text{max}}$}{
			Solve \eqref{eq: general discrete equation a.d.s. x dir.} to obtain $X^{1,n, \theta, \rho}$\;
			Normalize $X^{1,n, \theta, \rho}$ to satisfy \eqref{eq: X_norm = 1}\;
			Solve \eqref{eq: general discrete equation a.d.s. y dir.} to obtain $Y^{1,n, \theta, \rho}$\;
			\If{$\left( X^{1,n, \theta, \rho}, Y^{1,n, \theta, \rho} \right)$ satisfy \eqref{eq: criteria for tensor decomposition}}{
				break\;
			}
		}
		% $\overline{U^{n, \theta}}$ is solved in tensor decomposed form.\;
		Formulate \eqref{eq: discrete VMS fine: equ.} and \eqref{eq: discrete VMS fine: b.c.} into the form of \eqref{eq: general discrete equation a.d.s. x dir.} and \eqref{eq: general discrete equation a.d.s. y dir.}\;
		Initialize $X^{2,n, \theta, 0}_{q, J_x} \leftarrow \sin (q \pi J_x / N)$, $Y^{2,n, \theta, 0}_{q, J_y} \leftarrow \sin (q \pi J_y / N)$\;
		\For{$\rho = 1, \cdots, \rho_{\text{max}}$}{
			Solve \eqref{eq: general discrete equation a.d.s. x dir.} to obtain $X^{2,n, \theta, \rho}$\;
			Normalize $X^{2,n, \theta, \rho}$ to satisfy \eqref{eq: X_norm = 1}\;
			Solve \eqref{eq: general discrete equation a.d.s. y dir.} to obtain $Y^{2,n, \theta, \rho}$\;
			\If{$\left( X^{2,n, \theta, \rho}, Y^{2,n, \theta, \rho} \right)$ satisfy \eqref{eq: criteria for tensor decomposition}}{
				break\;
			}
		}
		Compute ${U^{n, \theta}}$ in its tensor-decomposed form.\;
		\If{$\overline{U^{n,\theta}}$ satisfies \eqref{eq: criteria for inter-scale iteration}}{
			break\;
		}
	}
	% $\overline{U^{n}}, U^n$ are solved in tensor decomposed form.\;
	\caption{VMS-TD Algorithm}
	\label{al: two-level VMS-TD}
\end{algorithm}
In numerical implementation, we take $\theta_{\text{max}} = 30$, $\rho_{\text{max}} = 50$.

\subsection{Extend to three levels} \label{sect: extend to three scales}

For the sake of brevity, let
\begin{equation}
    \overline{U^n} = \mathcal{C}(f, U_B^n, \bm{L}_1;\ U^n, \bm{I}_{12})
\end{equation}
represent the coarse-scale solution to \eqref{eq: discrete VMS coarse: equ.}\eqref{eq: discrete VMS coarse: b.c.}, where $f$ is the heat source function, $U_B^n$ is the boundary condition, $\bm{L}_1$ represents all coarse-scale grid information (including grid size, number of nodes, shape functions, etc.), $U^n$ is the fine-scale data, and $\bm{I}_{12}$ contains all the information connecting the coarse and fine scales (including the location of the fine-scale region, refinement ratio, etc.). Once all the parameters are given, \eqref{eq: discrete VMS coarse: equ.}\eqref{eq: discrete VMS coarse: b.c.} can be solved to obtain $\overline{U^n}$. This process is represented by $\mathcal{C}$. Similarly, the process of solving \eqref{eq: discrete VMS fine: equ.}\eqref{eq: discrete VMS fine: b.c.} is abbreviated as
\begin{equation}
    U^n = \mathcal{F}(f, \overline{U^n}, \bm{L}_2)
\end{equation}
where $f$ is still the heat source function, $\overline{U^n}$ provides the boundary condition, and $\bm{L}_2$ contains all the information related to fine-scale discretization.

% Thus, the two-level flowchart from $t_{n-1}$ to $t_n$ can be concisely represented by Algorithm \ref{al: two-level VMS}.

\begin{figure}[htbp]
	\centering
	\includegraphics[width=.6\textwidth]{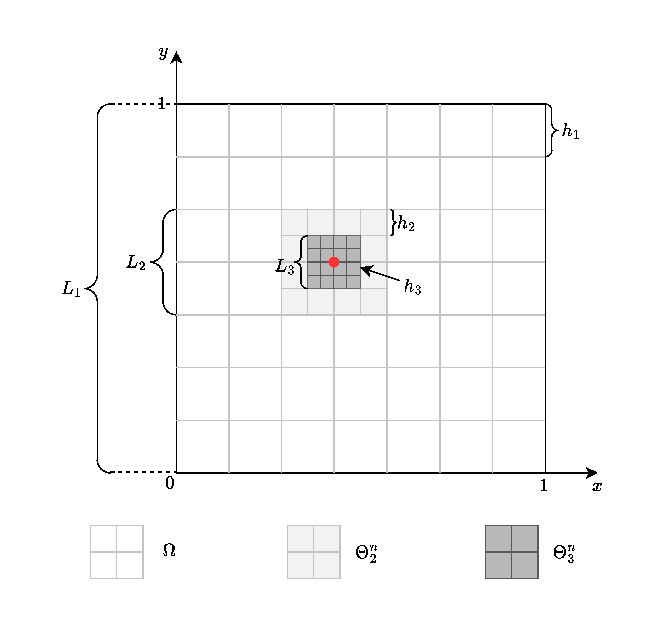}
	\caption{Schematical plot of a three-level mesh in 2D.}
	\label{fig: VMS 2D 3L mesh}
\end{figure}

For the three-level case, a region near the heat source is drawn as $\Theta_3^n \subset \Theta_2^n \subset \Omega$, as shown in Fig. \ref{fig: VMS 2D 3L mesh}. Let
\begin{equation}
	u^n = \overline{u^n} + (\widetilde{u^n})_2 + (\widetilde{u^n})_3 \quad \text{in} \quad \Omega,
\end{equation}
where
\begin{equation}
	(\widetilde{u^n})_2(x, y) = 0, \quad (x, y) \in \Omega \backslash \Theta_2^n;
\end{equation}
\begin{equation}
	(\widetilde{u^n})_3(x, y) = 0, \quad (x, y) \in \Omega \backslash \Theta_3^n.
\end{equation}

Define
\begin{equation}
	u_2^n = \overline{u^n} + (\widetilde{u^n})_2,
\end{equation}
and use $\overline{U^n}$, $U_2^n$ and $U_3^n$ to represent the discretized $\overline{u^n}$, $u_2^n$ and $u^n$, respectively. The three-level flowchart from $t_{n-1}$ to $t_n$ is shown in Algorithm \ref{al: three-level VMS}.

\begin{algorithm}[htbp]
	\SetAlgoLined
	
	Initialize $U_2^n \leftarrow 0$\;
	\While{$(U_1^{n}, U_2^n)$ have not converged}{
		$\overline{U^n} \leftarrow \mathcal{C}(f, U_B^n, \bm{L}_1;\ U_2^n, \bm{I}_{12})$\;
		Initialize $U_3^n \leftarrow 0$\;
		\While{$(U_2^{n}, U_3^n)$ have not converged}{
			$U_2^n \leftarrow \mathcal{C}(f, \overline{U^n}, \bm{L}_2;\ U_3^n, \bm{I}_{23})$\;
			$U_3^n \leftarrow \mathcal{F}(f, U_2^n, \bm{L}_3)$\;
		}
	}
	\caption{Three-level Algorithm}
	\label{al: three-level VMS}
\end{algorithm}

In practice, the iteration process can be accelerated by omitting the iteration between the second and third scales. This gives an accelerated three-level algorithm shown in Algorithm \ref{al: three-level VMS accelerated}.

\begin{algorithm}[htbp]
	\SetAlgoLined
	
	Initialize $U_2^n \leftarrow 0$ and $U_3^n \leftarrow 0$\;
	\While{$(U_1^{n}, U_2^n)$ have not converged}{
		$\overline{U^n} \leftarrow \mathcal{C}(f, U_B^n, \bm{L}_1;\ U_2^n, \bm{I}_{12})$\;
		$U_2^n \leftarrow \mathcal{C}(f, \overline{U^n}, \bm{L}_2;\ U_3^n, \bm{I}_{23})$\;
		$U_3^n \leftarrow \mathcal{F}(f, U_2^n, \bm{L}_3)$\;
	}
	\caption{Accelerated Three-level Algorithm}
	\label{al: three-level VMS accelerated}
\end{algorithm}
Both Algorithms \ref{al: three-level VMS} and \ref{al: three-level VMS accelerated} use the same stop criteria as the two-level case in Algorithm \ref{al: two-level VMS-TD}.

In most cases, the accelerated three-level algorithm converges. We use Algorithm \ref{al: three-level VMS accelerated} for all three-level numerical results presented in this paper.

\section{Numerical Results}\label{sect: numerical results}

In this section, we demonstrate the accuracy and efficiency of VMS-TD through several numerical examples.
First, we use a two-dimensional two-level example to demonstrate that VMS-TD reproduce faithfully TD with uniform fine mesh over the whole domain.
Then, a three-dimensional three-level example shows that the proposed algorithm applies to more levels and higher dimensional problems.
In addition, we compare the computational time between VMS-TD algorithm and TD with uniform fine mesh over the whole domain to show the significant efficiency improvement of the VMS-TD algorithm when dealing with large problems.

\subsection{Two-level VMS-TD simulations for 2D heat conduction}

\begin{figure}[h]
  \includegraphics[width=\textwidth]{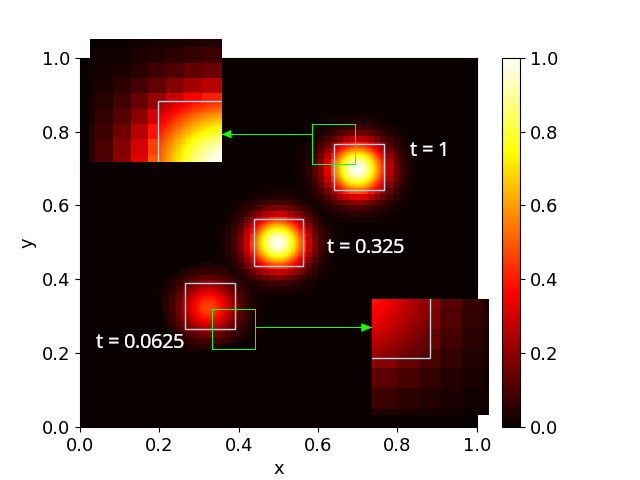}
	\caption{Two dimensional results $u(x,y,t)$ by two-level VMS-TD algorithm with $N_1 = N_2 = 64$, $Q=2$ at $t = 0.0625, 0.375, 1$.}
	\label{fig: 2D2LGenView}
\end{figure}

We choose $\Omega = [0, 1]^2$, $T = 1$, $\nu = 0.05$, and the heat source function
\begin{equation} \label{eq: 2D2L heat source}
	\begin{array}{rl}
	f(x, y, t) =&\bigg[\bigg( \dfrac{ (x - \mu_x(t)) \mu_x'(t)+(y - \mu_y(t)) \mu_y'(t) }{\sigma^2} -\nu\dfrac{(x-\mu_x(t))^2+(y-\mu_y(t))^2}{\sigma^4} +  \dfrac{2\nu}{\sigma^2} \bigg) \\[5mm]
&\qquad \left( 1 - e^{-\lambda t} \right) +\lambda e^{-\lambda t} \bigg]\exp\left( - \dfrac{(x - \mu_x(t))^2+(y - \mu_y(t))^2}{2\sigma^2} \right).
	\end{array}
\end{equation}
This gives an exact solution
\begin{equation}
	u_{\text{ex}}(x, y, t) = \exp\left( - \dfrac{(x - \mu_x(t))^2}{2\sigma^2} \right) \exp\left( - \dfrac{(y - \mu_y(t))^2}{2\sigma^2} \right) \left( 1 - e^{-\lambda t} \right),
\end{equation}
where $\mu_x(t) = \mu_y(t) = 0.3 + 0.4 t$, $\sigma= 0.05$, $\lambda = 10$.
Fig. \ref{fig: 2D2LGenView} shows the numerical results with $N_t = 128$, $N_1 = 64$, $N_2 = 64$, $L_2/L_1 = 0.125$. This gives $h_1=1/64,h_2=1/512$, and the coarsening ratio $h_1/h_2 = 8$. For each direction, we use two Gauss points for finite element analysis, and three Gauss points for error evaluation.

%\newpage
%\vspace{0.3in}
%\begin{figure}[h]
%\vspace{2.7in}
%	%\centering
%	\hspace{1in}
%  \includegraphics[width=20mm]{2D2LGenView_gimp.png}
%	\caption{Two dimensional results $u(x,y,t)$ by two-level VMS-TD algorithm with $N_1 = N_2 = 64$, $Q=2$ at $t = 0.0625, 0.375, 1$.}
%	\label{fig: 2D2LGenView}
%\end{figure}
%\vspace{0.3in}

For the time-marching scheme, we define an error
\begin{equation} \label{eq: def VMS-TD error 2D 2L}
	e_{\text{VMS-TD}} = \dfrac{ \dfrac{1}{N_t} \sum_{n = 1}^{N_t} \left\{ \int_{\Omega \backslash \Theta^n} \left[ \overline{u^n}(x, y) - u^n_{\text{ex}}(x, y) \right]^2 \ud (x, y) + \int_{\Theta^n} \left[ u^n(x, y) - u^n_{\text{ex}}(x, y) \right]^2 \ud (x, y) \right\}^{1/2} }{ \dfrac{1}{N_t} \sum_{n = 1}^{N_t}  \left\{ \dint_{\Omega}   \left[ u^n_{\text{ex}}(x, y) \right]^2 \ud (x, y) \right\}^{1/2} }.
\end{equation}
For comparison, we also compute by TD with uniform fine mesh over the whole domain. to solve \eqref{eq: time-discretized heat equ. G. weak form}, with fine grid covering the whole domain $\Omega$. That is to say, we use a uniform mesh size $h_2=\dfrac{L_2}{N_2}$ in each dimension for TD computation. Note that VMS-TD uses this mesh only for the subdomain and much coarse mesh size $h_1=\dfrac{L_1}{N_1}$ over the vast surrounding area. Therefore, we can not expect VMS-TD outperforming TD with uniform fine mesh over the whole domain. Moreover, difference between these two numerical solutions shows the influence of coarsening, i.e., the VMS strategy in our algorithm. For ease of presentation, we refer to this TD solution with uniform fine mesh over the whole domain as TD reference solution, and denote it by $u_{\text{ref}}$. The corresponding error is
\begin{equation} \label{eq: def TD error 2D 2L}
	e_{\text{ref}} = \dfrac{ \dfrac{1}{N_t} \sum_{n = 1}^{N_t} \left\{ \int_{\Omega} \left[ u_{\text{ref}}^n(x, y) - u^n_{\text{ex}}(x, y) \right]^2 \ud (x, y) \right\}^{1/2} }{ \dfrac{1}{N_t} \sum_{n = 1}^{N_t}  \left\{ \dint_{\Omega}   \left[ u^n_{\text{ex}}(x, y) \right]^2 \ud (x, y) \right\}^{1/2} }.
\end{equation}
Note that both errors are compared with the exact solution.

As is always for numerical simulations, there is a trade-off between accuracy and speed.
VMS-TD is faster than TD with uniform fine mesh over the whole domain, yet as discussed before, the accuracy quantified by $e_{\text{VMS-TD}}$ may not be better than $e_{\text{ref}}$ in general.

Now we explore how numerical parameters influence the differences between $e_{\text{VMS-TD}}$ and $e_{\text{ref}}$. In this way we elucidate the effect of the VMS formulation.

First, we fix the space discretization $N_1 = 64$ and $N_2 = 64$. The number of modes $Q = 2$ is also fixed. The time step size changes as $N_t = 4, 8, \cdots, 1024$.
\begin{figure}[h]
\setlength{\unitlength}{1mm}
\begin{picture}(120,2)\put(15,0){(a)}\put(100,0){(b)}\end{picture}
	\centering
\includegraphics[width=.5\textwidth]{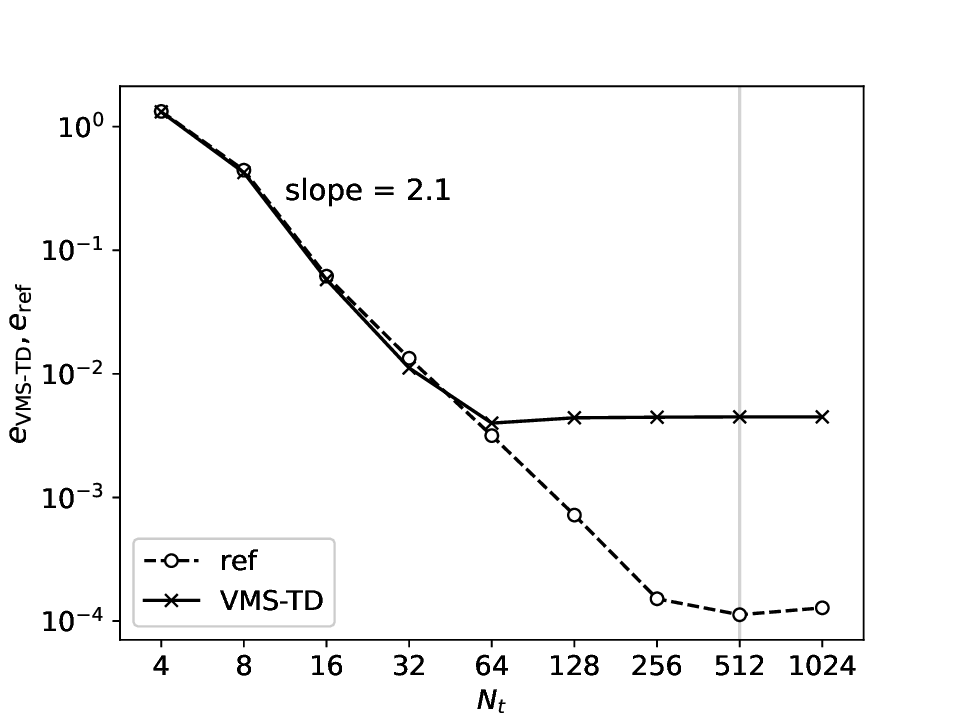}\includegraphics[width=.5\textwidth]{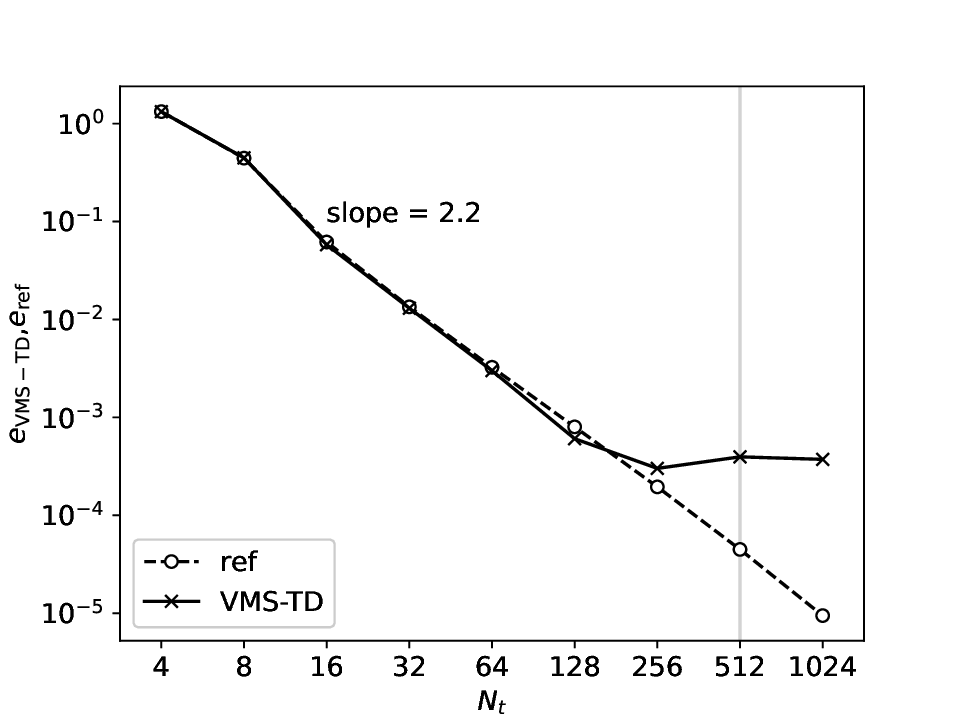}
	\caption{Numerical error $e_{\text{VMS-TD}}$ and $e_{\text{ref}}$ for $N_t = 4, 8, \cdots, 1024$, with $Q = 2$ and fixed mesh for two-dimensional test: (a)  $N_1 = N_2 = 64$; (b) $N_1=N_2=256$.}
	\label{fig: 2D2L L2Err-Nt}
\end{figure}
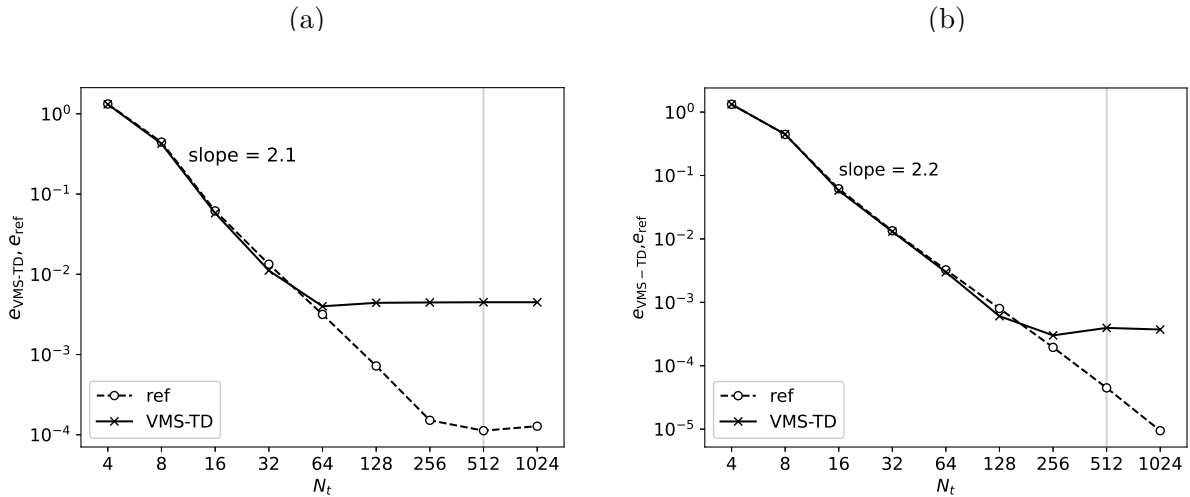

In the left subplot of Fig. \ref{fig: 2D2L L2Err-Nt}, we observe that for $N_t \leq 64$, VMS-TD faithfully reproduces TD with uniform fine mesh over the whole domain. This implies that the VMS coarsening strategy does not cause observable difference, and the temporal discretization error dominates. This is further verified by the convergence rate at about $2$, as the second-order time semi-discretization scheme is adopted in \eqref{eq: time-discretized heat equ.}.
For finer time integration, difference appears for VMS-TD versus TD with uniform fine mesh over the whole domain, where both errors are smaller than $10^{-2}$. Refining time integration does not improve the accuracy of VMS-TD for $N_t > 64$. In contrast, the TD reference solution with finer time step further improves for $N_t > 256$. This is because spatial discretization error dominates when time discretization is fine enough.

The accuracy may be improved by casting finer grid. In the right subplot of Fig. \ref{fig: 2D2L L2Err-Nt}, one order of magnitude improvement is reached by taking a mesh 4 times finer.

In the following explorations, we take $N_t = 512$ to guarantee that the time integration is accurate enough, thereby allows full consideration of spatial discretizations.

\begin{figure}[h]
	\centering
	\includegraphics[width=\textwidth]{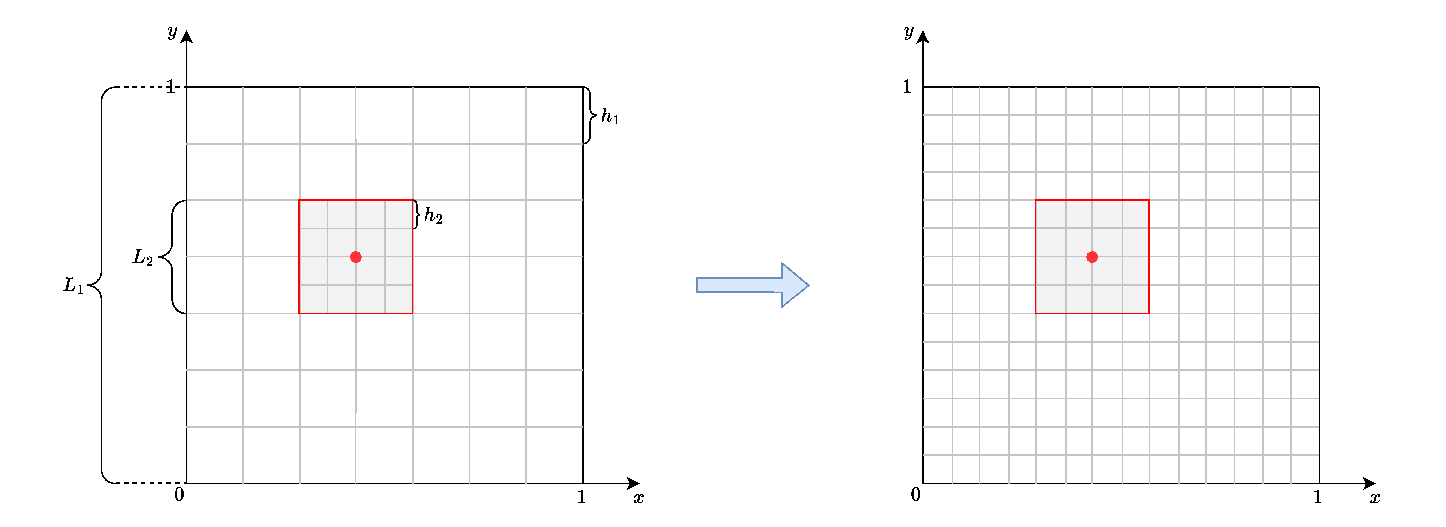}
	\caption{Refine the coarse-scale mesh.}
	\label{fig: refine_h1}
\end{figure}

Secondly, we fix the fine-scale mesh size $h_2$, while change the coarse scale mesh size $h_1$, as schematically shown in Fig. \ref{fig: refine_h1}.

\begin{figure}[h]
	\centering
	\includegraphics[width=.5\textwidth]{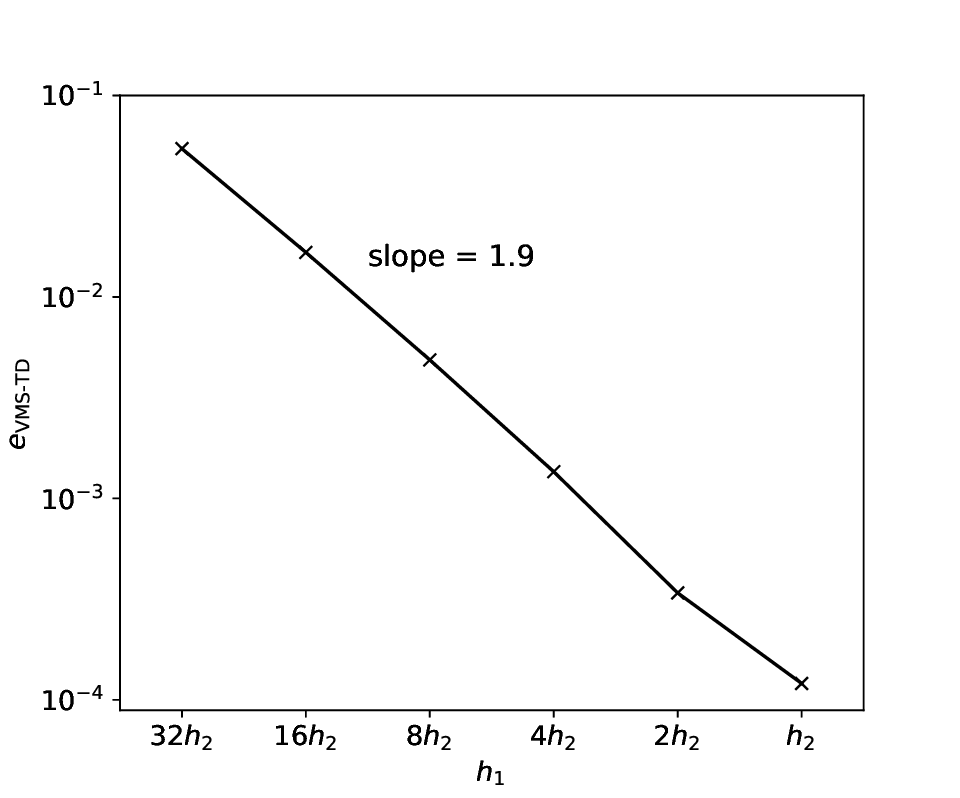}
	\caption{Numerical error $e_{\text{VMS-TD}}$ for different coarse mesh: $h_1 = 32 h_2, 16 h_2, \cdots, h_2$.}
	\label{fig: 2D2L L2Err-Nc}
\end{figure}

In Fig. \ref{fig: 2D2L L2Err-Nc}, when we refine the coarse grid size to the fine one, namely, $h_1=h_2$, VMS-TD reduces to TD with uniform fine mesh over the whole domain. Then the error is the same, at $1.5\times 10^{-4}$. On the other hand, for bigger $h_1$, the error increases with a slope $1.9$, due to the quadratic convergence of linear finite element analysis.

\begin{figure}[h]
	\centering
	\includegraphics[width=\textwidth]{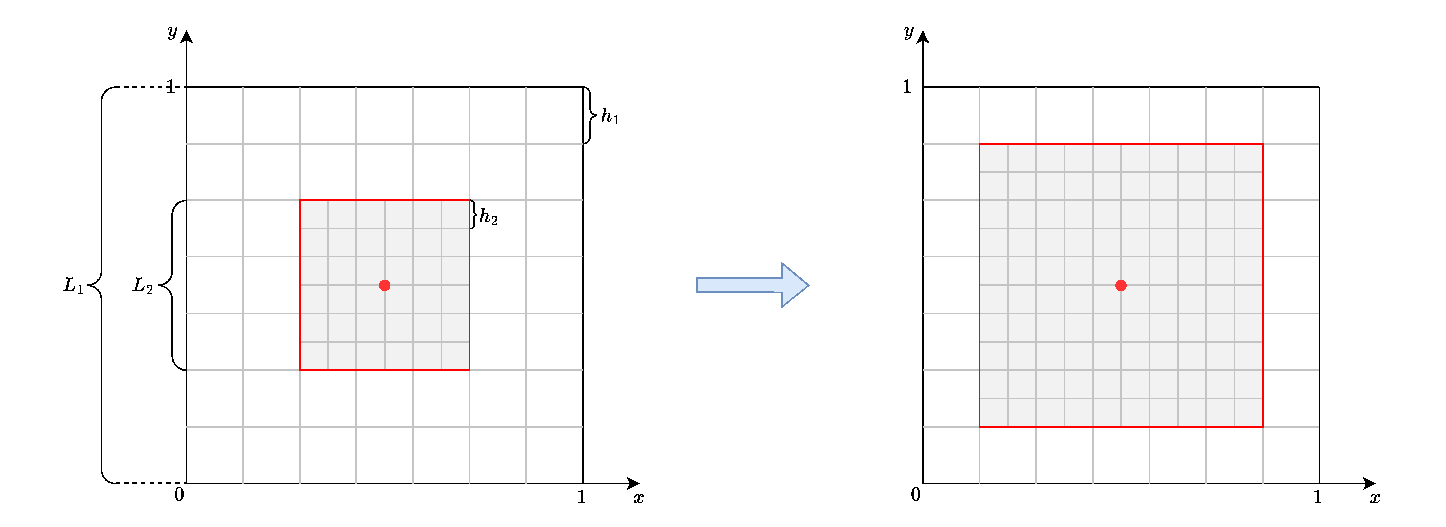}
	\caption{Enlarge fine-scale subdomain.}
	\label{fig: enlarge_L2.drawio}
\end{figure}

Thirdly, we fix $h_1=1/64$ and $h_2=1/512$, and enlarge the fine-scale subdomain $\Theta^n$ by changing $L_2$.
For this test, we fix the heat source centered at $\mu_x(t) = \mu_y(t) = 0.5$. See Fig. \ref{fig: enlarge_L2.drawio} for illustration.

\begin{figure}[h]
\setlength{\unitlength}{1mm}
\begin{picture}(120,2)\put(15,0){(a)}\put(100,0){(b)}\end{picture}
	\centering \includegraphics[width=.5\textwidth]{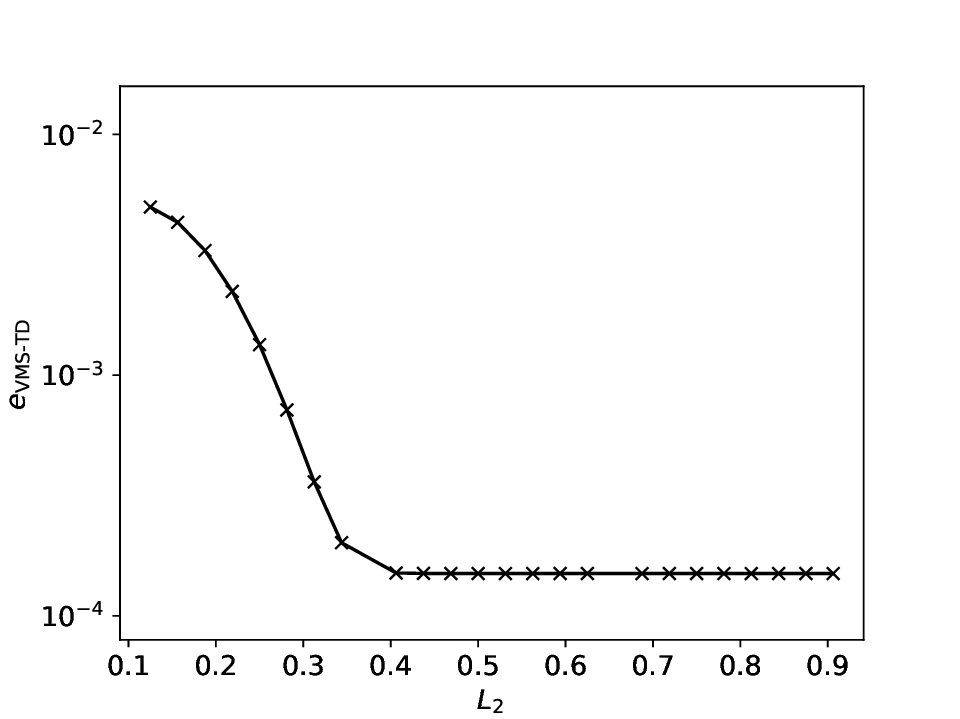}\includegraphics[width=.5\textwidth]{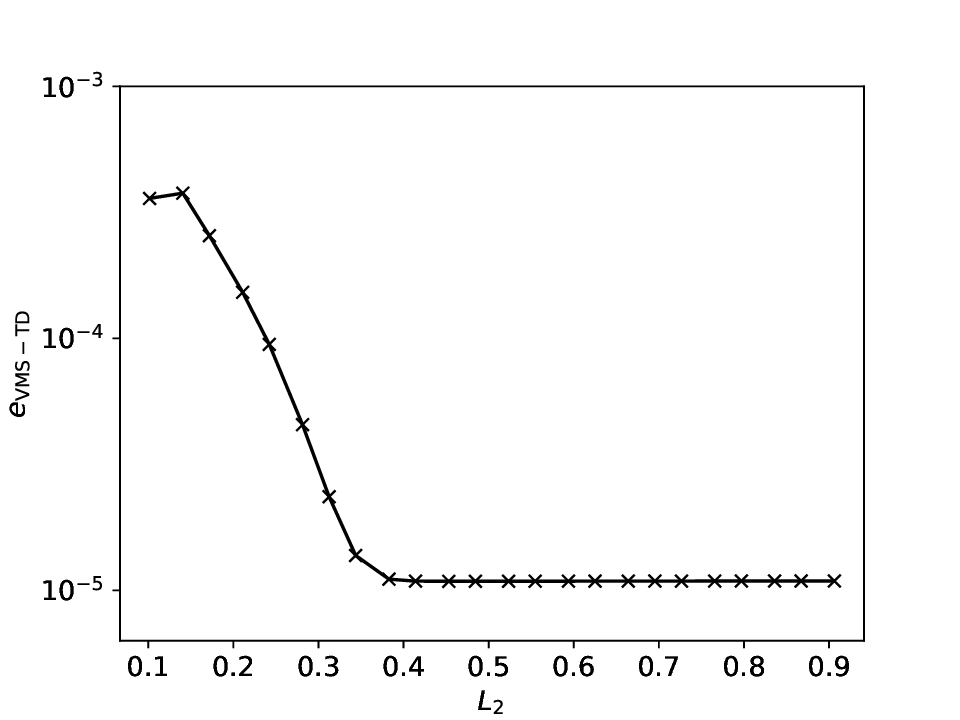}
	\caption{Numerical error $e_{\text{VMS-TD}}$ for $L_2 \in [0.1, 0.9]$, with fixed mesh size: (a) $h_1=1/64, h_2=1/512$; (b) $h_1=1/256,h_2=1/2048$.}
	\label{fig: 2D2L L2Err-lengthRatio12}
\end{figure}
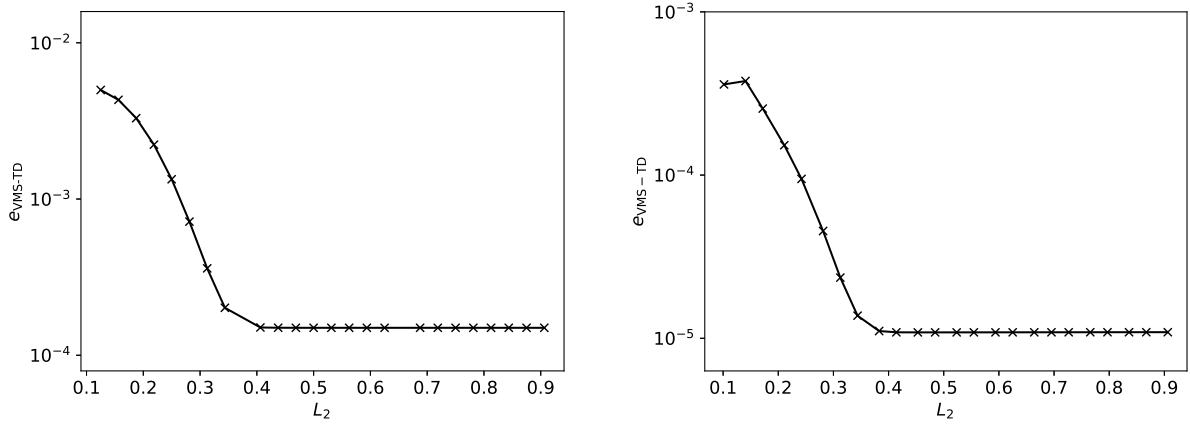

Numerical results in Fig. \ref{fig: 2D2L L2Err-lengthRatio12} shows that VMS-TD well reproduces the TD reference solution when $L_2 \geq 0.4 L_1 $, giving an error $e_{\text{VMS-TD}}=1.5\times 10^{-4}$ in the left subplot. This may be explained with the help of Fig. \ref{fig: 2D2LGenView}. Due to the heat conduction process, the solution is negligible outside of a box with side length $0.4$, centered at $(\mu_x, \mu_y)$. Therefore, when the fine scale TD computation correctly resolves solution in such a box, the numerical error in this box is small. Meanwhile, the temperature outside the box is small, so is the numerical error there. Consequently, the difference between the VMS-TD solution and the TD reference results over the entire computational domain is small. We remark that this threshold of $L_2$ is problem-specific, depending on the diffusivity coefficient and the form of heat source. Furthermore, the right subplot with finer grid improves the accuracy by one order of magnitude. Because of the TD formulation, this is reached with only 4 times many DoFs, in contrast to $4^2=16$ times DoFs for direct numerical simulation with finite element analysis.

Results for refining the coarse mesh and enlarging the fine subdomain are quantitatively summarized in Table \ref{tab: L2Err for different h1 and L2}. Choices that yield error smaller than $0.1\%$ are bold-faced.
When $h_1$ is small enough or $L_2$ is big enough, the accuracy is rather high. For moderate $h_1$ and $L_2$, this table suggests suitable combined choices.

\newpage

\vspace{10mm}
\setlength{\unitlength}{1mm}
\linethickness{1pt}
\begin{picture}(50,-5)
\put(35,0){\vector(1,0){100}}
\put(5,-5){\vector(0,-1){30}}
\end{picture}
\vspace{-10mm}
\begin{table}[H]
	\centering

	\begin{tabular}{l|cccccc}
	\hline
	$\quad \quad h_1$ & $32h_2$ & $16h_2$ & $8h_2$ & $4h_2$ & $2h_2$ & $h_2$ \\
	$L_2 \quad$ &  &   &  &   &  &  \\
	\hline
	0.125 & $8.512\%$ & $2.027\%$ & $0.499\%$ & $0.132\%$ & $0.039\%$ & $\bm{0.015\%}$ \\
	\hline
	0.188 &  $1.531\%$ & $1.302\%$ & $0.431\%$ & $0.111\%$ & $\bm{0.032\%}$ & $\bm{0.015\%}$ \\
	\hline
	0.25 & $1.531\%$ & $0.501\%$ & $0.223\%$ & $\bm{0.071\%}$ & $\bm{0.023\%}$ & $\bm{0.015\%}$ \\
	\hline
	0.3 & $\bm{0.071\%}$ & $0.117\%$ & $\bm{0.036\%}$ & $\bm{0.017\%}$ & $\bm{0.015\%}$ & $\bm{0.015\%}$ \\
	\hline
	0.35 & $\bm{0.071\%}$ & $\bm{0.023\%}$ & $\bm{0.016\%}$ & $\bm{0.015\%}$ & $\bm{0.015\%}$ & $\bm{0.015\%}$ \\
	\hline
	0.4 & $\bm{0.018\%}$ & $\bm{0.015\%}$ & $\bm{0.015\%}$ & $\bm{0.015\%}$ & $\bm{0.015\%}$ & $\bm{0.015\%}$ \\
	\hline
	\end{tabular}
	\caption{Numerical error at different coarse mesh size $h_1$ and fine mesh subdomain size $L_2$. Refining coarse mesh (indicated by the horizontal arrow), and enlarge fine mesh subdomain (indicated by the vertical arrow) reduce the numerical error.}
	\label{tab: L2Err for different h1 and L2}
\end{table}

\subsection{Three-level VMS-TD simulations for 3D heat conduction}

We present three-level VMS-TD simulation results for heat conduction in 3D to further illustrate the efficiency of the proposed algorithm. Choose $\Omega = [0, 1]^3$, $T = 1$, $\nu = 0.05$, and take a heat source function
\begin{equation}
	\begin{array}{rl}
	f(x, y, t) = & \bigg[ \bigg(  \dfrac{ (y - \mu_y(t)) \mu_y'(t) }{\sigma^2} - \nu\dfrac{(x - \mu)^2+(y - \mu_y(t))^2+(z-\mu)^2 }{\sigma^4}+ \dfrac{3\nu}{\sigma^2} \bigg) \\[5mm] & \qquad \left(1 - e^{-\lambda t} \right)+ \lambda e^{-\lambda t} \bigg] \exp\left( - \dfrac{(x - \mu)^2+(y - \mu_y(t))^2+(z-\mu)^2 }{2\sigma^2} \right),
	\end{array}
\end{equation}
with $\mu =0.5$, $\mu_y(t) = 0.3 + 0.4 t$, $\sigma = 0.02$, $\lambda = 10$. The heat source moves along $y$ direction.

The exact solution is
\begin{equation}
	u_{\text{ex}}(x, y, t) =  \exp\left( - \dfrac{(x - \mu)^2+(y - \mu_y(t))^2+(z-\mu)^2 }{2\sigma^2} \right)\left(1-e^{\lambda t}\right).
\end{equation}

\begin{figure}[h]
%\vspace{3.5in}
%	%\centering
	\includegraphics[width=\textwidth]{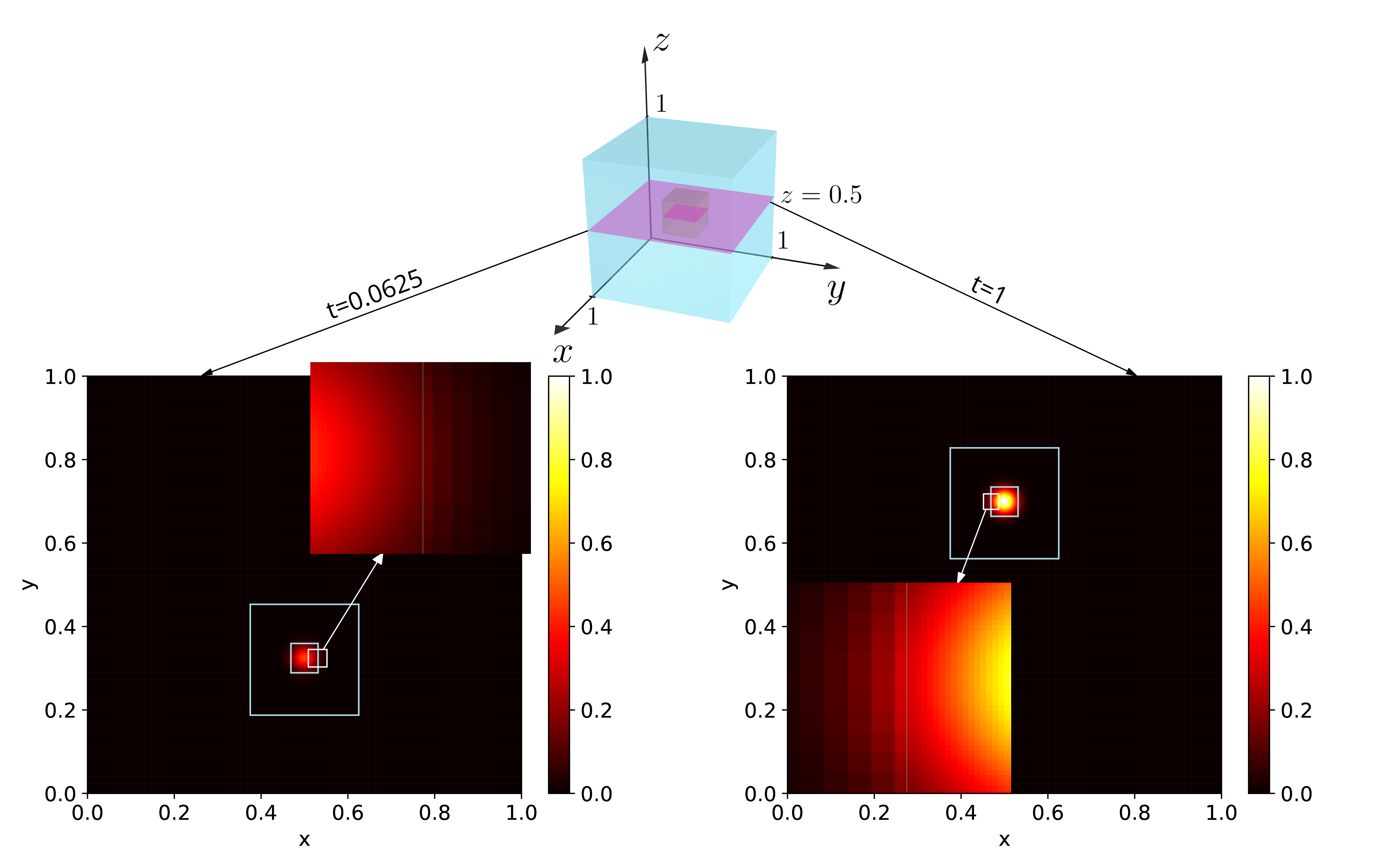}
	\caption{Three dimensional results $u(x,y,0.5,t)$ by three-level VMS-TD algorithm with $N_1 = N_2 = N_3 = 64$, $Q=2$ at different time: (a) $t=0.0625$; (b) $t=0.1$.}
	\label{fig: 3D3LGenView}
\end{figure}

Fig. \ref{fig: 3D3LGenView} shows the results with $N_t = 512$, $N_1 = N_2 = N_3 = 64$, $L_2= 1/4, L_3 = 1/16$. The corresponding mesh sizes are $h_1=1/64,h_2 =1/256,h_3=1/1024$, and coarsening ratios are $h_1/h_2=h_2/h_3=4$. For each direction, we use two Gauss points for finite element analysis, and three Gauss points for error evaluation.
Note that we only plot results on the plane $z = 0.5$ in Fig. \ref{fig: 3D3LGenView}.

Errors for the multi-level VMS-TD and  with uniform fine mesh over the whole domain are defined similar to the two-level case in \eqref{eq: def VMS-TD error 2D 2L} and \eqref{eq: def TD error 2D 2L}.

First, we fix the space discretization and number of modes, with $N_1 =N_2 = N_3 = 64$ and $Q = 2$. The time step size changes with $N_t = 4, 8, \cdots, 2048$.

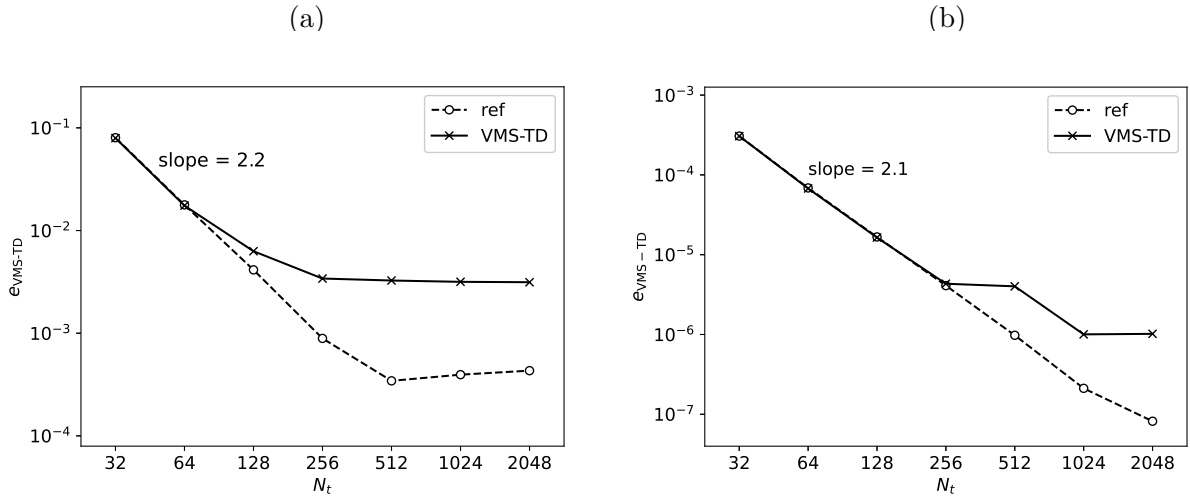
\begin{figure}[H]
\setlength{\unitlength}{1mm}
\begin{picture}(120,2)\put(15,0){(a)}\put(100,0){(b)}\end{picture}
	\centering	\includegraphics[width=.5\textwidth]{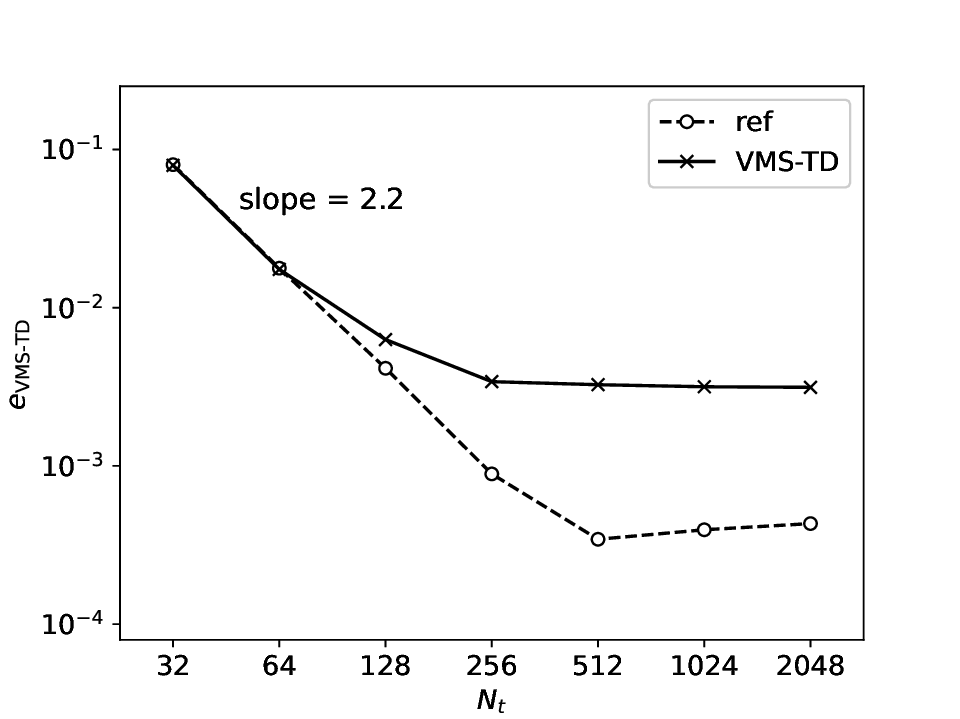}\includegraphics[width=.5\textwidth]{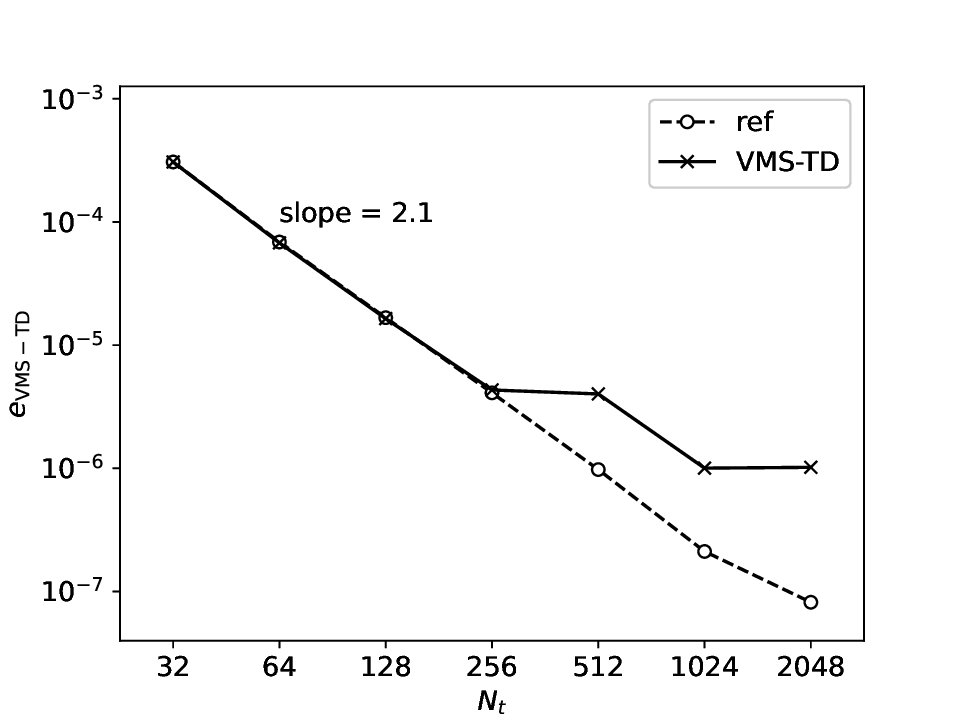}
	\caption{ Numerical error $e_{\text{VMS-TD}}, e_{\text{ref}}$ for $N_t = 32, 64, \cdots, 2048$,with $Q = 2$ and fixed mesh for three dimensional test: (a) $N_1 = N_2 =N_3= 64$; (b) $N_1=N_2=N_3=256$.}
	\label{fig: 3D3L L2Err-Nt}
\end{figure}

VMS-TD offers a special coarsening of TD with uniform fine mesh over the whole domain, hence $e_{\text{VMS-TD}}\geq e_{\text{ref}}$ in general. In the left subplot of Fig. \ref{fig: 3D3L L2Err-Nt}, when $N_t \leq 64$, the temporal error dominates, hence the coarsening in VMS-TD does not cause observable differences, and the convergence rate is about $2$ for the second-order time semi-discretization scheme \eqref{eq: time-discretized heat equ.}. On the other hand, for finer time integration, difference appears. As a matter of fact, for the time step size no smaller than $1/{N_t}=1/64$, refining spatial mesh does not improve the accuracy, because the time integration dominates numerical error. When we refine the time integration, the spatial discretization error now dominates, causing relative error at the level of $10^{-2}$ for VMS-TD. In contrast, refining time integration for TD with uniform fine mesh over the whole domain further improves accuracy before $N_t = 512$.  Of course TD with uniform fine mesh over the whole domain is numerically more expensive, which will be discussed in the next subsection separately. We further notice that a four times finer grid leads to more than three orders of magnitude accuracy improvement, as illustrated by the right subplot of Fig. \ref{fig: 3D3L L2Err-Nt}. This demonstrates that the accuracy of VMS-TD improves along with grid refinements.

In all following explorations, we take $N_t = 512$ to guarantee that time integration is fine enough, allowing full consideration of spatial discretizations.

\begin{figure}[H]
	\centering
	\includegraphics[width=.5\textwidth]{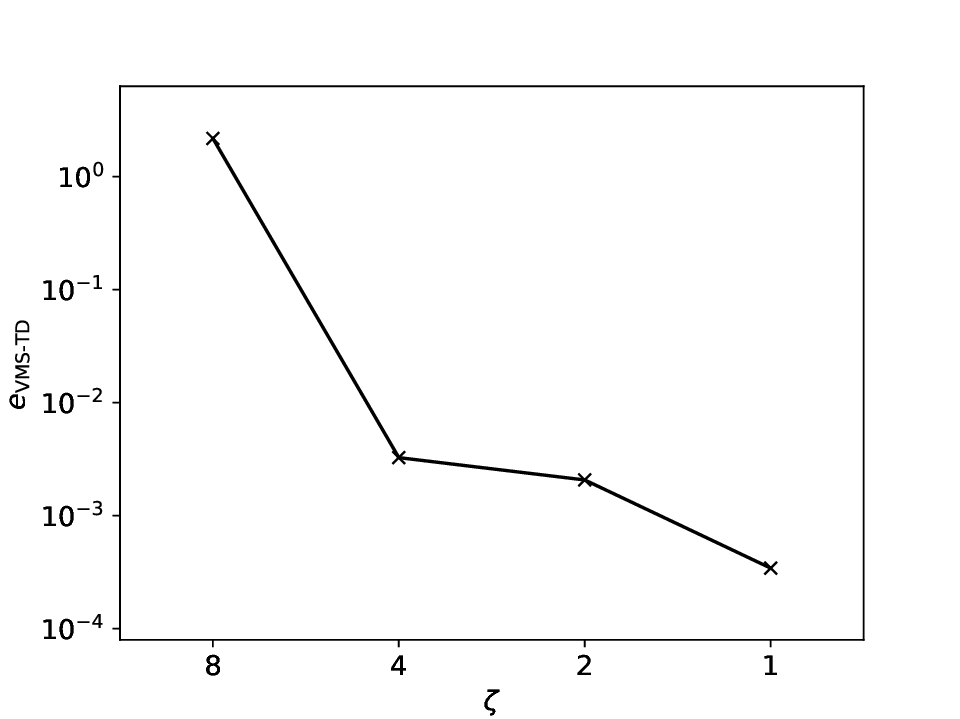}
	\caption{Three dimensional results: error $e_{\text{VMS-TD}}$ for coarsening ratio $\zeta = 1,2,4,8$.}
	\label{fig: 3D3L L2Err-Nc}
\end{figure}

Secondly, we fix the level-3 mesh size $h_3$, while change the level-1 and level-2 mesh sizes $h_1$, $h_2$ with coarsening ratio $\zeta$, i..e., $ h_1=\zeta h_2, h_2=\zeta h_3$. The relative error $e_{\text{VMS-TD}}$ is plotted in Fig. \ref{fig: 3D3L L2Err-Nc}.
When the level-1 mesh size $h_1$ and level-2 mesh size $h_2$ reaches the level-3 mesh size $h_3$, VMS-TD and TD with uniform fine mesh over the whole domain are the same. Therefore we observe that the error reaches $3.16\times 10^{-4}$ at $\zeta = 1$, the same as that for the TD reference solution (cf. left subplot of Fig. \ref{fig: 3D3L L2Err-Nt} for $N_t=512$). Though we only obtain results until $\zeta = 8$, it is clear that the three-level results for 3D heat equation has the same trend as the two-level result for 2D heat equation.

Thirdly, we fix $h_1=1/64,h_2 =1/256,h_3=1/1024$, and enlarge the level-2 and level-3 subdomains $\Theta_2^n$ and $\Theta_3^n$ by changing $L_2$ and $L_3$ simultaneously (keep $L_3 = L_2^2$). For this test, we fix the heat source centered with $ \mu_y(t) = 0.5$.

\begin{figure}[H]
\setlength{\unitlength}{1mm}
\begin{picture}(120,2)\put(15,0){(a)}\put(100,0){(b)}\end{picture}
	\centering	\includegraphics[width=.5\textwidth]{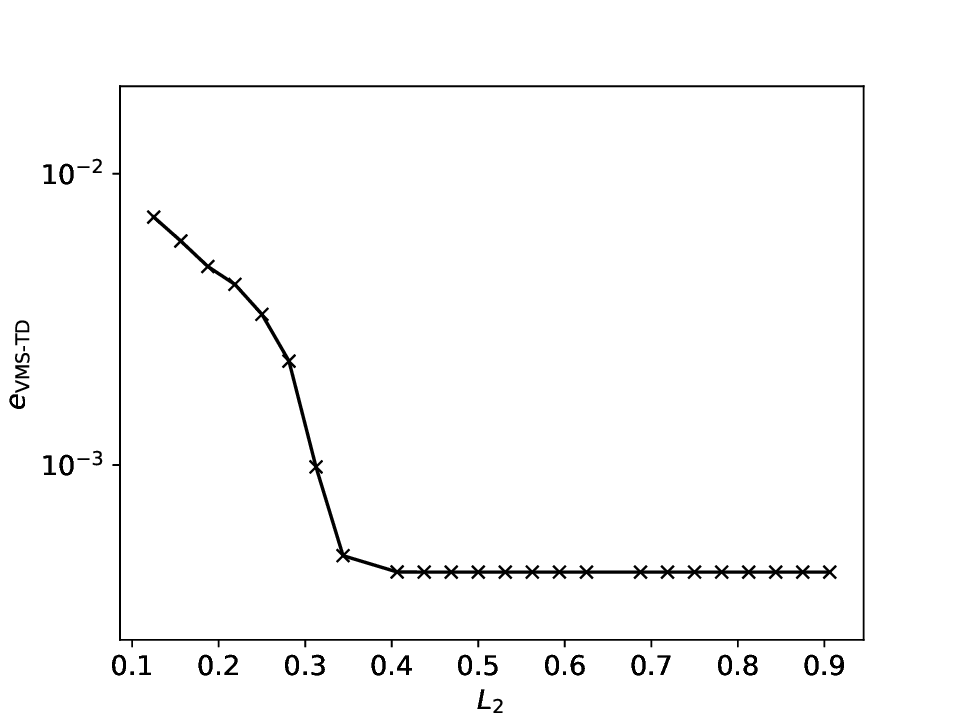}\includegraphics[width=.5\textwidth]{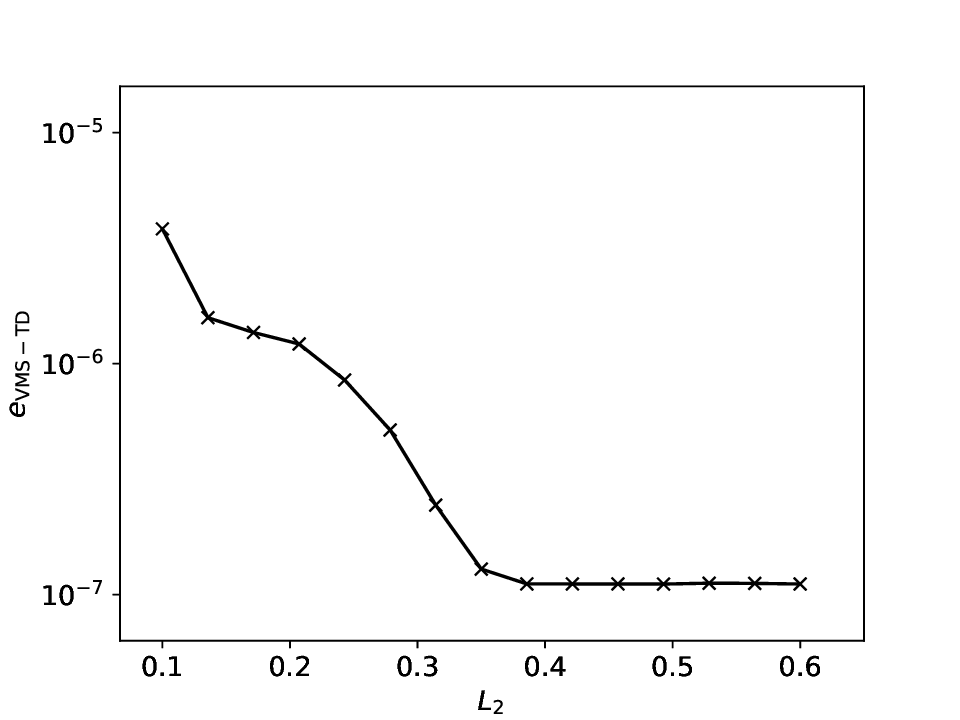}
	\caption{Numerical error $e_{\text{VMS-TD}}$ with fixed mesh for three dimensional test: (a) $N_1=N_2=N_3=64,h_1=1/64,h_2=1/256,h_3=1/1024$ for $L_2 \in [0.1, 0.9]$; (b) $N_1=N_2=N_3=256,h_1=1/256,h_2=1/1024,h_3=1/4096$ for $L_2 \in [0.1, 0.6]$.}
	\label{fig: 3D3L L2Err-lengthRatio12}
\end{figure}
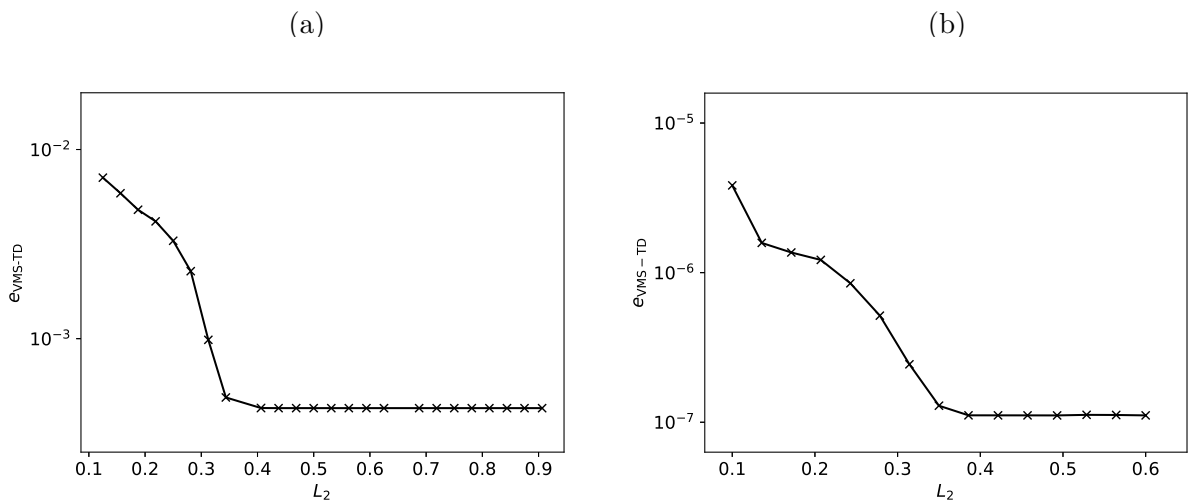

Numerical results in Fig. \ref{fig: 3D3L L2Err-lengthRatio12} shows that VMS-TD reproduces full TD when $L_2 \leq 0.4$ (when $L_3 \leq 0.16$). This may be explained with Fig. \ref{fig: 3D3LGenView}, where the solution is negligible outside of a cube with side length $0.16$, centered at $(\mu, \mu_y(t), \mu)$. When the TD reference solution well approximates the exact solution in such a cube, so does the VMS-TD solution. We remark that this threshold of $L_2$ is problem-specific, depending on the diffusivity coefficient and the form of heat source. We further notice that a four times finer grid leads to more than three orders of magnitude accuracy improvement, as illustrated by the right subplot of Fig. \ref{fig: 3D3L L2Err-lengthRatio12}.
Again we see that Fig. \ref{fig: 3D3L L2Err-lengthRatio12} has the same behavior as Fig. \ref{fig: 2D2L L2Err-lengthRatio12}, demonstrating the effectiveness of extension for the two-level approach to three-level.

\subsection{Efficiency} \label{subsec: Efficiency}
The merit of the multi-level VMS-TD algorithm lies in the efficiency as compared to the TD with uniform fine mesh over the whole domain. While the TD algorithm is already advantageous in reducing DoFs and numerical costs, we compare three-level VMS-TD results for the 3D heat conduction problem to demonstrate its power and potential.

Since VMS-TD is a time-marching method, we focus on one time step to evaluate the complexity and computational cost. For the sake of simplicity, we change
the level-1 grid number $N_1$ (in each dimension), and take $N_2 = N_3 = N_1$, $h_1 = \dfrac{1}{N_1}$, $h_2 = \dfrac{10}{N_1^2}$, and $h_3 = \dfrac{10^2}{N_1^3}$. This gives $L_1=1, L_2=N_2h_2=\dfrac{10}{N_1}, L_3=N_3h_3=\left(\dfrac{10}{N_1}\right)^2$. We note that the TD with uniform fine mesh over the whole domain reaches the same resolution  using $N_{\text{ref}} = \dfrac{N_1^3}{10^2}$ nodes in each dimension.

Furthermore, if taking direct numerical simulation with finite difference/finite element analysis, the DoFs reach $N_{\text{ref}}^3 = \dfrac{N_1^9}{10^6}$. Even if an elaborated full multi-grid method is adopted, the computational cost is then $O\left(\dfrac{N_1^9}{10^6}\right)$, and the convergence rate reaches that of the underlying finite difference method. For heat conduction with moving heat source in LPBF, the advantage of the proposed VMS-TD is clear. Unlike multi-grid method which is purely a technique for solving algebraic equations, here we take into account the physical process and pursue fine-scale resolution only for a very small subdomain where it is necessary.

\begin{figure}[H]
	\centering
	\includegraphics[width=\textwidth]{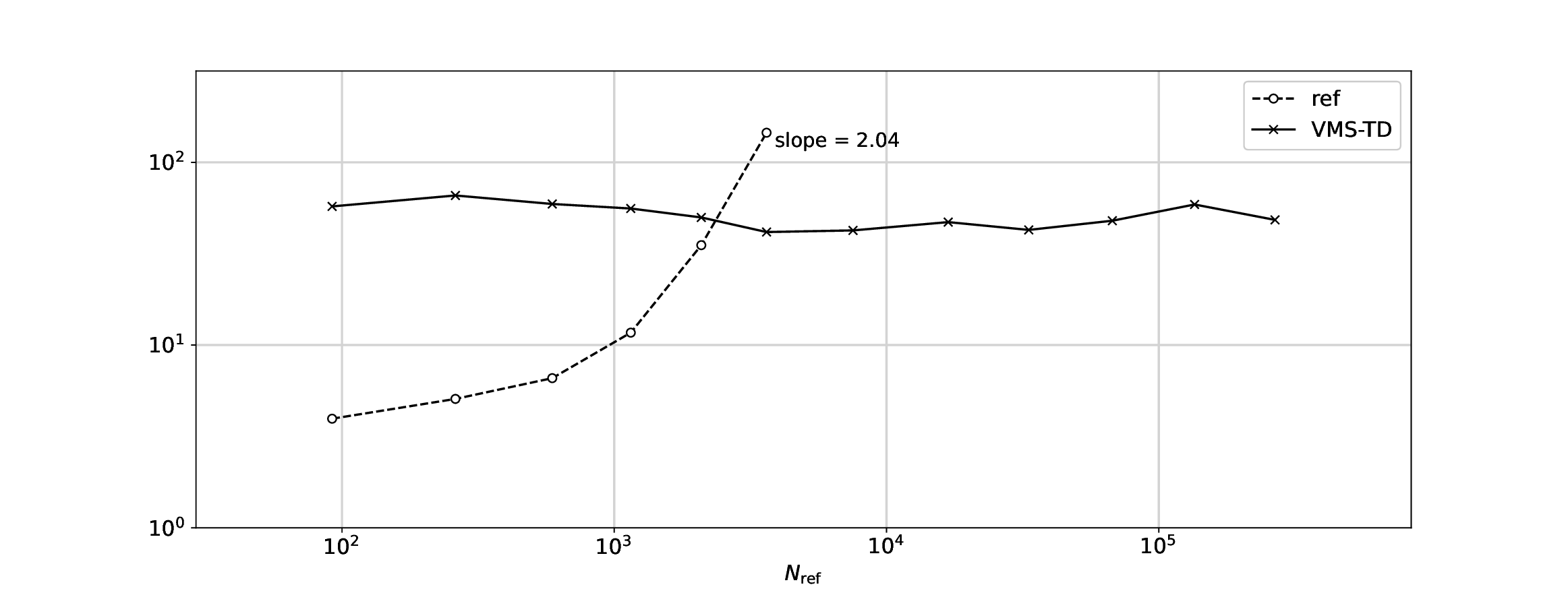}
	\caption{Average wall clock time (seconds) for a single time step with respect to different $N_{\text{ref}}$ in each dimension.}
	\label{fig: lgWCT-lgN_2ull}
\end{figure}

Fig. \ref{fig: lgWCT-lgN_2ull} shows the average wall clock time for computing each time step.

All numerical simulations are run on the same GPU (Nvidia GeForce RTX 4060 Ti). As we can see in Fig. \ref{fig: lgWCT-lgN_2ull}, TD with uniform fine mesh over the whole domain is faster for small $N_{TD}$ (big mesh size), while three-level TD outperforms for large $N_{TD}$ (small mesh size).

In  TD with uniform fine mesh over the whole domain, the major numerical costs come from solving the linear system, with complexity $O(N_{\text{ref}}^3)$.
% In our implementation, we use the CUDA library 'cusolver' to reach the rate $2.04$, which is smaller than $3$.
In the three-level VMS-TD approach, we solve linear systems of size $N_1=N_2=N_3  \sim O( N_{\text{ref}}^{1/3} )$, each with complexity  complexity  $O\left((N_{\text{ref}}^{1/3})^3 \right) = O(N_{\text{ref}})$.
The number of such linear systems to solve depends on the number of iterations, including iterations in TD and those across the scales.
In Fig. \ref{fig: lgWCT-lgN_2ull}, wall clock time of the three-level VMS-TD approach remains virtually unchanged for different meshes. This is an exciting result, indicating its potential in large engineering simulations.

% Because of the limitation of GPU memory, we cannot give out full-mesh computation for very large $N_{\text{ref}}$, but we can estimate by extending the dashed-line. If we extend the dashed-line using the slope $2.04$, it takes $10^{6} \text{s} \approx 278 \text{h}$ to compute a single time step for $N_{\text{ref}} = 267300$. Using three-level approach, it takes only $50$s, which is about $2 \times 10^4$ times faster. The acceleration will be more obvious for even larger number of meshes.

\begin{figure}[H]
	\centering
	\subfigure[]{
	\includegraphics[width=0.48\textwidth]{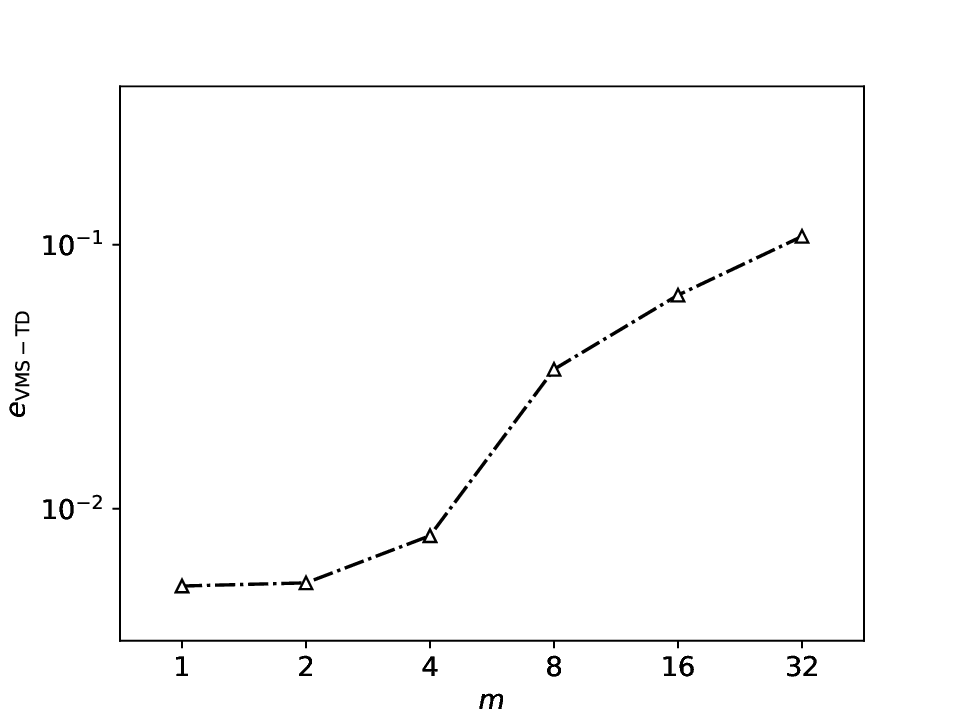}
	\label{fig: 2D2L L2Err-lgTSR(a)}
	}
	\subfigure[]{
	\includegraphics[width=0.48\textwidth]{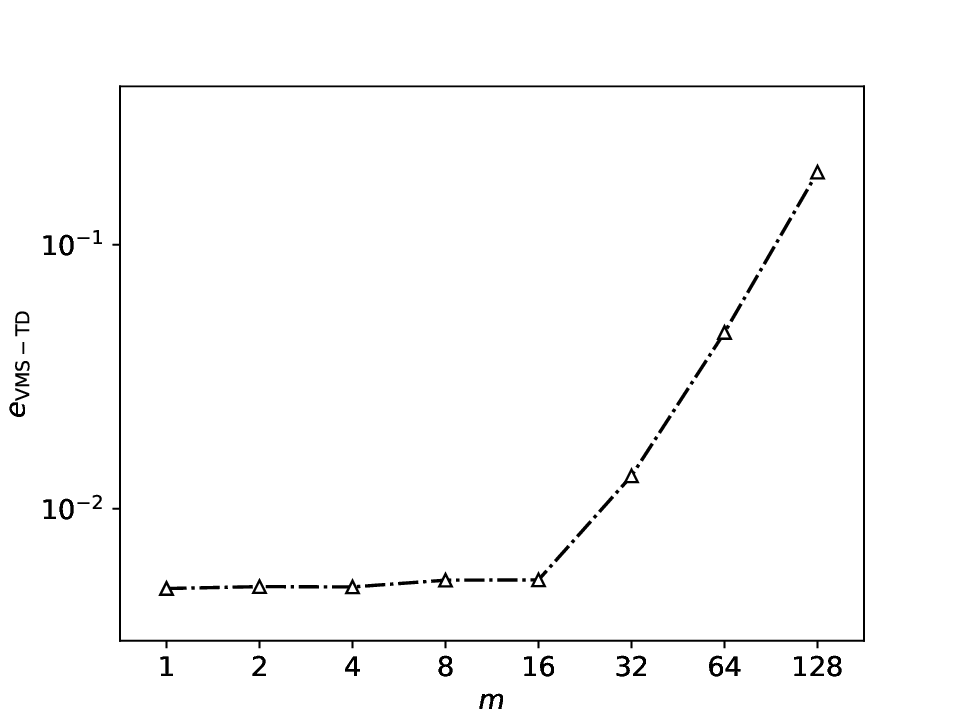}
	\label{fig: 2D2L L2Err-lgTSR(b)}
	}
	\caption{Numerical error $e_{\text{VMS-TD}}$ for time step coarsening ratio $m$: (a) $N_t = 128$; (b) $N_t = 512$}
	\label{fig: 2D2L L2Err-lgTSR}
\end{figure}

Now we make two remarks. First, our computations are based on dense matrix, and the largest problem we can solve on a single GPU (4060 Ti) is a $4096 \times 4096 \times 4096$ space mesh by TD with uniform fine mesh over the whole domain. In contrast, the three-level VMS-TD algorithm can solve a problem corresponding to a $1,000,000 \times 1,000,000 \times 1,000,000$ uniform fine mesh. If using sparse matrix, this can be further enlarged. Second, the VMS-TD algorithm can also incorporate the multiple time step technique. See Appendix for details. We show some preliminary results for the 2D heat conduction with a fixed heat source at $ \mu_y(t) = 0.5$, where all other parameters are kept the same as used for Fig. \ref{fig: 2D2LGenView}. While the fine scale computations are performed with $\Delta t=\dfrac{T}{N_t}$, the coarse scale computations are performed with time step size $m\Delta t$. See Fig. \ref{fig: 2D2L L2Err-lgTSR} for $N_t = 128$ and 512. The above VMS-TD corresponds to the case $m=1$. We observe that for relatively small $m$, the accuracy is preserved despite of the coarse time integration. Effective coupling of multiple time step integration and multi-level VMS-TD algorithm is an on-going study, and shall be presented in future work.

\section{Conclusion} \label{sect: conclusion}
In this paper, we propose a multi-level VMS-TD algorithm by combining the VMS formulation and TD representation to solve effectively each time step in a second-order time marching scheme. One key challenge is to avoid data-reconstruction and keep all data in arrays. To reach high efficiency in implementation, we put forward special treatments for interface and moving subdomain. This pieces up the two scales seamlessly. More precisely, for the interface, we split it into lines and corners, where a variable-separated function can be defined. For the moving subdomain, we split the updated data on the irregular shaped subdomain into several rectangles, supporting variable-separated forms. Numerical results show the effectiveness of the proposed algorithm. Moreover,it is extended to three-level and higher space dimensions. Numerical examples illustrate its efficiency and accuracy. Multiple time step integration is also incorporated into the algorithm. We demonstrate that multi-level VMS-TD algorithm greatly reduces computational cost, and has the potential to reach a very fine resolution for engineering problems.

Our approach still has many limitations. Basically, it is based on a time-marching scheme, where we only decompose the space into variable-separated form. It then takes a lot of time steps to solve the heat equation. It is interesting and challenging to combine VMS and TD based on a space-time method, which solves the whole equation in one stroke and avoids accumulative errors. In addition, TD requires the computational domain at each scale to be rectangular/cubic. How to deal with domain of arbitrary shape is a task demanding further explorations.

\section*{Funding} This research did not receive any specific grant from funding agencies in the public, commercial, or not-for-profit sectors.
%
%\section*{Data Availability} Enquiries about data availability should be directed to the authors.
%
%\section*{Declarations}
%\textbf{Conflict of interest} The authors have not disclosed any competing interests.
%

\section*{Appendix: Multiple time step integration}

\begin{figure}[htbp]
	\centering
	\includegraphics[width=.6\textwidth]{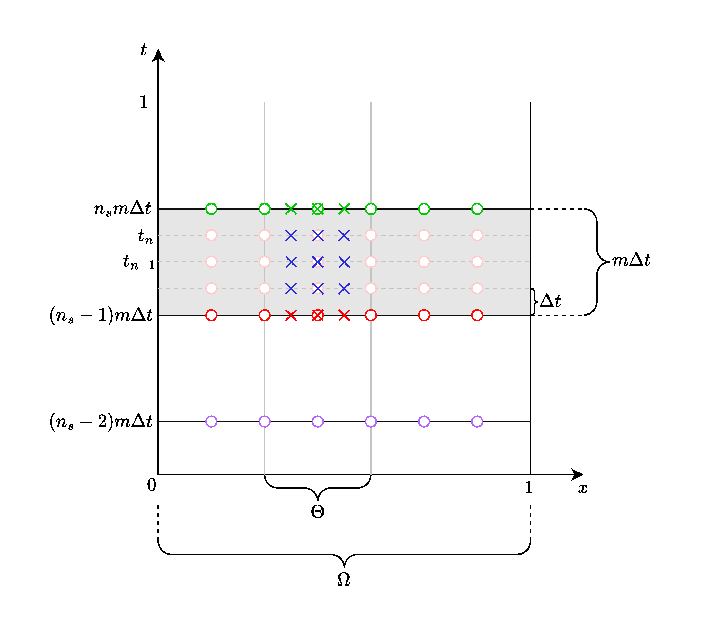}
	\caption{Multiple time steps}
	\label{fig: multiple time steps}
\end{figure}

We may incorporate the multiple time step integration scheme to the multi-level VMS-TD algorithm. See Fig. \ref{fig: multiple time steps} for illustration.

Choosing a time step ratio $m$, we group $m$ steps into a slab (grey area). We now describe the time integration over one slab. With temperature at the coarse and fine scales given at the beginning of the slab (red circles and crosses), we calculate temperature at the slab end  (green circles and crosses). First, we take all green circles the same as red circles as an initial guess, and use a quadratic interpolation (using green, red, purple circles) to get pink circles. Then we only solve the fine-scale equations step by step to get blue crosses, noting that in this process the coarse values (pink circles) remain unchanged. After getting the final step fine-scale values (green crosses), we solve the coarse-scale equation to couple the two scales, in turn the green circles are updated. We then iterate within this slab until the results converge.

The algorithm in the $n_s$-th slab is shown in Algorithm \ref{al: multiple time step}.
\begin{algorithm}[htbp]
	\SetAlgoLined
	
	Initialize $\overline{U^{n_s m} }  \leftarrow  \overline{U^{(n_s-1) m} }$ \;
	\While{$\overline{U^{n_s m} }$ have not converged}{
		Get $\overline{U^{(n_s-1) m + 1}}, \cdots, \overline{U^{n_s m - 1} }$ by quadratic interpolation \;
		\For{$n = (n_s-1) m + 1, \cdots, n_s m$}{
			$U^n \leftarrow \mathcal{F}(f, \overline{U^n}, \bm{L}_2)$\;
		}
		$\overline{U^{n_s m}} \leftarrow \mathcal{C}(f, U_B^{n_s m}, \bm{L}_1;\ U^{n_s m}, \bm{I}_{12})$
	}
	\caption{Multiple time step VMS-TD Algorithm}
	\label{al: multiple time step}
\end{algorithm}

Noting that this algorithm use starting values of the previous slab (purple circles), we conserve the original time-marching scheme to compute values in the first slab. We use quadratic interpolation because our method is based on a second-order time-marching scheme.

\end{document}